\documentclass[12pt]{article}
\usepackage{amssymb}
\usepackage{amsmath}

\input{tcilatex}

\begin{document}

\begin{center}
{\LARGE Monge-Kantorovich norms}

{\LARGE on spaces of vector measures}

\bigskip

Ion Chi\c{t}escu $^{a}$, Radu Miculescu $^{b,\ast }$, Lucian Ni\c{t}\u{a} $%
^{c}$, Loredana Ioana $^{d}$

\bigskip

$^{a}$ \textit{Faculty of Mathematics and Computer Science, University of
Bucharest, Academiei Str. 14, 010014, Bucharest, Romania}

$^{b}$ \textit{Faculty of Mathematics and Computer Science, University of
Bucharest, Academiei Str. 14, 010014, Bucharest, Romania}

$^{c}$ \textit{Technical University of Civil Engineering, Lacul Tei Blvd.,
122-124, 020396, Bucharest, Romania}

$^{d}$ \textit{Faculty of Mathematics and Computer Science, University of
Bucharest, Academiei Str. 14, 010014, Bucharest,}

\textit{\ Romania}
\end{center}

\medskip

------------------------------------------------------------------------------------------------

\textbf{Abstract}

\medskip

One considers Hilbert space valued measures on the Borel sets of a compact
metric space. A natural numerical valued integral of vector valued
continuous functions with respect to vector valued functions is defined.
Using this integral, different norms (we called them Monge-Kantorovich norm,
modified Monge-Kantorovich norm and Hanin norm) on the space of measures are
introduced, generalizing the theory of (weak) convergence for probability
measures on metric spaces. These norms introduce new (equivalent) metrics on
the initial compact metric space.

\medskip



\medskip

---------------------------------------------------------------------------------------------

* Corresponding author

\textit{E-mail addresses}: ionchitescu@yahoo.com (Ion Chi\c{t}escu),
miculesc@yahoo.com (Radu Miculescu), luci6691@yahoo.com (Lucian Ni\c{t}\u{a}%
), loredana.madalina.ioana@gmail.com (Loredana Ioana)\newpage

\textbf{1. Introduction}

\bigskip

We introduce and study several metrics on certain spaces of vector measures.
This is done using a vector integral (with numerical values) previously
introduced by us.

A short history of the problem follows. The story began long ago, with the
problem of mass transport initiated by G. Monge in 1781 (how to fill up a
hole with the material from a given pile of sand, in an optimal way, i.e.
with a minimal cost, see [18]). The problem was, actually, very difficult
and G. Monge proposed a complicated geometrical solution. Many years after,
in 1887, P. Appell completely solved the problem with complicated
variational methods (see [1]). L. V. Kantorovich, inventor of linear
programming, attacked the problem in a totally different way. First he
considered the discreet variant of the problem, "embedding" it in the theory
of linear programming and totally solving it (see [12]) in a way suitable
for successful use of computers (these results constitute a major part of
the reasons for the Nobel prize - of course for economy - received by L.V.
Kantorovich). Afterwards, he transformed the problem in an abstract way,
working for a compact metric space instead of a finite set and for a measure
instead of a vector in $\mathbb{R}^{n}$. Alone or jointly with his student
G.S. Rubinstein (see [14] and [15]), he succeeded in completely solving the
new, abstract problem. The necessary mathematical tools were the theory of
normed spaces and different metrics on spaces of measures, let us call them
Kantorovich-Rubinstein metrics. In the treatises [13] and [21] these facts
are clearly explained, with many details.

The study of different metrics on spaces of measures in closely related to
the theory of convergence of probability measures (especially weak
convergence) on \ metric spaces, where the Lipschitz functions play an
important role (see [2], [6] and [19]). It is within the framework of this
theory that the formalism of Kantorovich-Rubinstein-type metrics appears
more clearly.

Our main goal in the present paper is to extend the theories of metrics and
convergence in spaces of probabilities (or scalar measures) to (similar)
theories in spaces of vector measures. The most suitable framework seemed to
us to be the framework of Hilbert space valued measures. So, let $X$ be a
Hilbert space.

We needed first an integral. Consequently, we elaborated the theory of a
numerical valued integral of continuous functions on a compact metric space
taking values in $X$, with respect to a measure of bounded variation taking
values in $X$. Our integral is sesquilinear (not bilinear in the complex
case) and uniform. In the second part of the paragraph dedicated to
preliminaries we expose (without proofs) the main properties of this
(natural) integral, among them being some computing devices and an
antilinear (not linear in the complex case) and isometric isomorphism
between the space of measures and the dual of the space of continuous
functions. Details and proofs will appear in "Sesquilinear uniform vector
integral", Proceedings - Mathematical Sciences.

Having this integral, we pass in the main paragraph ("Results") of the paper
to the study of the space of $X$-valued measures (of bounded variation)
defined on the Borel sets of a compact metric space $(T,d)$. Different
metrics (and locally convex topologies) are introduced on these spaces, a
major role being played by the Lipschitz functions (the use of these
functions makes the new introduced metrics to be "topologically sensitive",
as one can see in the last subparagraph).

We begin with the Monge-Kantorovich norm and the corresponding metric. We
continue with the weak$^{\ast }$ topology. It is seen that, in case of
finite dimensional $X$, the weak$^{\ast }$ and the Monge-Kantorovich
topologies coincide, which is not the case for infinite dimensional $X$, as
one can see later. Afterwards, we introduce a new norm on a subspace of the
space of $X$-valued measures and the corresponding metric on some subsets of
the whole space of $X$-valued measures. We called them "the modified
Monge-Kantorovich norm ",\ respectively "the modified Monge-Kantorovich
metric". The modified Monge-Kantorovich norm and Monge-Kantorovich norm are
equivalent. The results described up to now constitute generalizations (for
measures) of many results concerning probability measures.

The next subparagraph generalizes (for vector measures) the ideas of L.
Hanin (see [11] and [10]), the ideas being to "extend" the modified
Monge-Kantorovich norm to the whole space of $X$-valued measures (actually
one obtains a norm which is equivalent to the modified Monge-Kantorovich
norm on the initial subspace). Using the previous results, we introduce a
counterexample showing that weak$^{\ast }$ topology and Monge-Kantorovich
topology are not the same for infinite dimensional spaces.

The last subparagraph has a somewhat different character. Namely, using the
previously introduced metrics on the space of measures, we can equip the
underlying compact metric space $T$ with the corresponding new metrics and
we show that all these new metrics are equivalent to the initial metric $d$
(which is not the case of the variational metric).

We feel obliged to add that the use of "Monge-Kantorovich" name which we
preferred seemed suitable to us for historical reasons. Many times the
"Kantorovich-Rubinstein" name would have been more correct.

The authors hope to use the results in the present paper in a subsequent
paper dedicated to applications (e.g. fractals).

\bigskip

\textbf{2. Preliminary facts}

\bigskip

\textbf{Notations and general notions}

\bigskip

We begin with some notations and general notions appearing throughout the
paper.

As usual $\mathbb{N}=\{1,2,...,n,...\}=$ the non null positive integers, $K=$
the scalar field (either $K=\mathbb{R}$ or $K=\mathbb{C}$), $\mathbb{R}%
_{+}=\{x\in \mathbb{R\mid }x\geq 0\}$, $\overline{\mathbb{R}_{+}}=\mathbb{R}%
_{+}\cup \{\infty \}$, $K^{n}=\{(x_{1},x_{2},...,x_{n})\mid x_{i}\in K\}$,
where $n\in \mathbb{N}$.

For an arbitrary set $T$, we write $\mathcal{P}(T)=\{A\mid A\subset T\}$ and
for $A\subset T$, $\varphi _{A}:T\rightarrow K$ will be the characteristic
(indicator) function of $A$ acting via $\varphi _{A}(t)=0$ if $t\notin A$
and $\varphi _{A}(t)=1$ if $t\in A$.\ If $A\subset T$, the complementary set
of $A$ is $C_{D}=\{t\in T\mid t\notin A\}$. We write $(a_{i})_{i\in
I}\subset T$ (or $(a_{i})_{i}\subset T$) to denote the fact that the family $%
(a_{i})_{i\in I}$ has the property $a_{i}\in T$ for any $i\in I$. In
particular, one can consider sequences $(a_{n})_{n\in \mathbb{N}}\subset T$
(or $(a_{n})_{n}\subset T$). In this case, when we write $%
(a_{n_{p}})_{p}\subset (a_{n})_{n}$, this means that $(a_{n_{p}})_{p}$ is a
subsequence of $(a_{n})_{n}$.

If $f:X\rightarrow Y$ is a function, we shall often write $x\mapsto f(x)$ to
designate the fact that the image of $x\in X$ under $f$ is $f(x)$. Assuming $%
f$ is injective, the (generalized) inverse of $f$ is the function $%
f^{-1}:f(X)\rightarrow X$ acting via $f^{-1}(y)\overset{def}{=}x$, where $%
x\in X$ is uniquely determined by the condition $f(x)=y$. Let $X,Y,Z$ be
three sets and $A\subset X$, $B\subset Y$. Let $f:A\rightarrow Y$ be such
that $f(A)\subset B$ and let $g:B\rightarrow Y$. The (generalized)
composition of $f$ and $g$ is the function $g\circ f:A\rightarrow Z$ acting
via $(g\circ f)(x)=g(f(x))$.

Let $X$ be a vector space over $K$. For any $x\in X$, we write $Sp(x)$ for
the vector space generated by $x$, i.e. for $\{\alpha x\mid \alpha \in K\}$.
If $f:T\rightarrow K$ a function, we can define the function $%
fx:T\rightarrow K$ acting via $fx(t)=f(t)x$.

Now, let us consider a topological space $(T,\tau )$. If $(a_{n})_{n\in 
\mathbb{N}}\subset T$ and $a\in T$, we write $a_{n}\underset{n}{\rightarrow }%
a$ to designate the fact that the sequence $(a_{n})_{n}$ converges to $a$.
Supplementarily, let us consider (more generally) a preordered set $(\Delta
,\leq )$ ($u\leq u$, for any $u\in \Delta $; $u\leq v$ and $v\leq w$ implies 
$u\leq w$, for any $u,v,w\in \Delta $) which is directed (for any $u$, $v$
in $\Delta $ there exists $w\in \Delta $ such that $u\leq w$ and $v\leq w$).
We consider a function $f:\Delta \rightarrow T$, write $f(\delta )=x_{\delta
}$ for any $\delta \in \Delta $ and identify $f\equiv (x_{\delta })_{\delta
\in \Delta }$. Under these circumstances, we write $(x_{\delta })_{\delta
\in \Delta }$ net $T$ (or $(x_{\delta })_{\delta }$ net $T$). Let also $x\in
T$. Then we write $x_{\delta }\underset{\delta }{\rightarrow }x$ (and we say
that $(x_{\delta })_{\delta }$ converges to $x$) if for any (basic)
neighborhood $V$ of $x$ there exists $\delta (V)\in \Delta $ such that $%
x_{\delta }\in V$ whenever $\delta \in \Delta $, $\delta \geq \delta (V)$.
For any $A\subset T$ and any $a\in T$, we have the equivalence: $a\in 
\overline{A}$ (the closure of $A$) if and only if there exists $(a_{\delta
})_{\delta \in \Delta }$ net $A$ such that $a_{\delta }\underset{\delta }{%
\rightarrow }a$.

If $(T,d)$ is a metric space, $\emptyset \neq A\subset T$ and $x\in T$, the
distance from $x$ to $A$ is defined via $d_{1}(x,A)=\inf \{d(x,a)\mid a\in
A\}$ (clearly $d_{1}(x,\{a\})=d(x,a)$). Then $d_{1}(x,A)=0$ if and only if $%
x\in \overline{A}$. For $x$ and $y$ in $T$ one has $\left\vert
d_{1}(x,A)-d_{1}(y,A)\right\vert \leq d(x,y)$. Considering two non empty
sets $A\subset T$ and $B\subset T$, the distance between $A$ and $B$ is
defined via $\delta (A,B)=\inf \{d(a,b)\mid a\in A,b\in B\}$. Clearly $%
\delta (\{a\},B)=d_{1}(a,B)$. If $\emptyset \neq A\subset T$, the diameter
of $A$ is $diam(A)=\sup \{d(x,y)\mid x,y\in A\}$.

For any normed space $(X,\left\Vert .\right\Vert )$, the dual of $X$ is $%
X^{^{\prime }}=\{V:X\rightarrow K\mid V$ is linear and continuous$\}$. Then $%
X^{^{\prime }}$ becomes a Banach space, when equipped with the (operator)
norm $\left\Vert V\right\Vert _{0}=\sup \{\left\vert V(x)\right\vert \mid
x\in X,\left\Vert x\right\Vert \leq 1\}$. Usually we write only $X$ (instead
of $(X,\left\Vert .\right\Vert )$) in order to designate a normed space.

If $X$ is a Hilbert space, we shall write $(x\mid y)$ for the scalar product
of the elements $x,y\in X$. Hence the scalar product $(.\mid .)$ yields the
norm $\left\Vert .\right\Vert $, acting via $\left\Vert x\right\Vert =\sqrt{%
(x\mid x)}$.

The space $K^{n}$ becomes (canonically) a Hilbert space with the scalar
product $(x\mid y)=\underset{i=1}{\overset{n}{\sum }}x_{i}\overline{y_{i}}$,
where $x=(x_{1},x_{2},...,x_{n})$ and $y=(y_{1},y_{2},...,y_{n})$. Hence $%
\left\Vert x\right\Vert \geq \left\vert x_{i}\right\vert $ for any $i\in
\{1,2,...,n\}$, where $x=(x_{1},x_{2},...,x_{n})$. The space of sequences $%
l^{2}=\{x=(x_{n})_{n\in \mathbb{N}}\mid x_{n}\in K$, $\overset{\infty }{%
\underset{n=1}{\sum }}\left\vert x_{n}^{2}\right\vert <\infty \}$ becomes a
Hilbert space with the scalar product $(x\mid y)=\underset{n=1}{\overset{%
\infty }{\sum }}x_{n}\overline{y_{n}}$, where $x=(x_{n})_{n\in \mathbb{N}}$
and $y=(y_{n})_{n\in \mathbb{N}}$. Let us note that $\left\Vert x\right\Vert
=(\overset{\infty }{\underset{n=1}{\sum }}\left\vert x_{n}^{2}\right\vert )^{%
\frac{1}{2}}$, where $x=(x_{n})_{n\in \mathbb{N}}$.

For general topology, see [16] and [7]. For functional analysis, see [8],
[13] and [20].

\bigskip

\textbf{A sesquilinear uniform integral}

\bigskip

In this subparagraph, we shall present without proofs the sesquilinear
uniform integral which will be used throughout the paper.

Let $(T,d)$ be a compact metric space and $X$ a Hilbert space with scalar
product $(.\mid .)$ and corresponding norm $\left\Vert x\right\Vert =\sqrt{%
(x\mid x)}$. The Borel sets of $T$ will be $\mathcal{B}\subset \mathcal{P}%
(T) $. The vector space $C(X)=\{f:T\rightarrow X\mid f$ is continuous$\}$ is
a Banach space with norm $f\mapsto \left\Vert f\right\Vert =\sup
\{\left\Vert f(t)\right\Vert \mid t\in T\}$. Actually, $C(X)$ is a closed
space of the Banach space $B(X)=\{f:T\rightarrow X\mid f$ is bounded$\}$
equipped with the norm $f\mapsto \left\Vert f\right\Vert =\sup \{\left\Vert
f(t)\right\Vert \mid t\in T\}$ (the confusional same notation $\left\Vert
f\right\Vert $ for $f\in C(X)$ and $f\in B(X)$ is justified).

A function $f:T\rightarrow X$ will be called \textit{simple} if it has the
form $f=\overset{m}{\underset{i=1}{\sum }}\varphi _{A_{i}}x_{i}$, with $%
A_{i}\in \mathcal{B}$ and $x_{i}\in X$ (one can always consider that the
sets $A_{i}$ are mutually disjoint and $\underset{i=1}{\overset{m}{\cup }}%
A_{i}=T$). One has $S(X)=\{f:T\rightarrow X\mid f$ is simple$\}\subset B(X)$
and $S(X)$ is a vector subspace.

A function $\mu :\mathcal{B}\rightarrow X$ is called a $\sigma $-additive
measure if $\mu (\underset{n=1}{\overset{\infty }{\cup }}A_{n})=\underset{n=1%
}{\overset{\infty }{\sum }}\mu (A_{n})$ for any sequence $(A_{n})_{n}\subset 
\mathcal{B}$ of mutually disjoint sets. For such a $\mu $ and $A\in \mathcal{%
B}$, one can define \textit{the variation of }$\mu $\textit{\ over }$A$ as
follows. We shall say that a finite family $(A_{i})_{i\in
\{1,2,...,m\}}\subset \mathcal{B}$ is \textit{a partition of }$A$ if the
sets $A_{i}$ are mutually disjoint and $\underset{i=1}{\overset{m}{\cup }}%
A_{i}=A$. Then the variation of $\mu $ over $A$, denoted by $\left\vert \mu
\right\vert (A)$, is defined via $\left\vert \mu \right\vert (A)=\sup \{%
\underset{i=1}{\overset{m}{\sum }}\left\Vert \mu (A_{i})\right\Vert \mid
(A_{i})_{i\in \{1,2,...,m\}}$ is a partition of $A\}$ (the supremum is taken
for all possible partitions of $A$). We say that $\mu $\textit{\ is of
bounded variation }if\textit{\ }$\left\vert \mu \right\vert (T)<\infty $.

The vector space (with natural operations) $cabv(X)=\{\mu :\mathcal{B}%
\rightarrow X\mid \mu $ is a $\sigma $-additive measure of bounded variation$%
\}$ becomes a Banach space when equipped with the norm $\mu \mapsto
\left\Vert \mu \right\Vert \overset{def}{=}\left\vert \mu \right\vert (T)$.
Any $\sigma $-additive measure $\mu :\mathcal{B}\rightarrow K$ is of bounded
variation. The topology on $cabv(X)$ generated by this norm will be called
the variational topology and will be denoted by $\mathcal{T}(var,X)$. For
any $a\in (0,\infty )$, let $B_{a}(X)=\{\mu \in cabv(X)\mid \left\Vert \mu
\right\Vert \leq a\}$. Then $\mathcal{T}(var,X)$ induces the topology $%
\mathcal{T}(var,X,a)$ on $B_{a}(X)$. For a sequence $(\mu _{n})_{n}\subset
cabv(X)$ and for $\mu \in cabv(X)$, $\mu _{n}\overset{\text{var}}{\underset{n%
}{\rightarrow }}\mu $ means that $(\mu _{n})_{n}$ converges to $\mu $ in $%
\mathcal{T}(var,X)$. Notice that $\mu _{n}\overset{\text{var}}{\underset{n}{%
\rightarrow }}\mu $ $\Rightarrow $ $\mu _{n}\overset{\text{u}}{\underset{n}{%
\rightarrow }}\mu $, the last symbol denoting uniform convergence.
Consequently, $\mu _{n}\overset{\text{var}}{\underset{n}{\rightarrow }}\mu $
implies $\mu _{n}(A)\underset{n}{\rightarrow }\mu (A)$ for any $A\in 
\mathcal{B}$. Notice that, if $\mu \in cabv(K)$ and $x\in X$, then $\mu x\in
cabv(X)$ and $\left\Vert \mu x\right\Vert =\left\Vert \mu \right\Vert
\left\Vert x\right\Vert $.

In the same way, if $(f_{n})_{n}\subset B(X)$ and $f\in B(X)$, we write $%
f_{n}\overset{\text{u}}{\underset{n}{\rightarrow }}f$ to denote the fact
that $(f_{n})_{n}$ converges uniformly to $f$ (i.e. $(f_{n})_{n}$ converges
to $f$ in the Banach space $B(X)$). The closure of $S(X)$ in $B(X)$ is the
space of \textit{totally measurable functions} denoted by $TM(X)$. So $TM(X)%
\overset{def}{=}\overline{S(X)}$. One has $C(X)\subset TM(X)$.

Now, let $\mu \in cabv(X)$. For any $f=\overset{m}{\underset{i=1}{\sum }}%
\varphi _{A_{i}}x_{i}\in S(X)$, \textit{the integral of }$f$\textit{\ with
respect to }$\mu $ is $\int fd\mu \overset{def}{=}\overset{m}{\underset{i=1}{%
\sum }}(x_{i}\mid \mu (A_{i}))$ (the definition does not depend upon the
representation of $f$). Because $\left\vert \int fd\mu \right\vert \leq
\left\Vert \mu \right\Vert \left\Vert f\right\Vert $, the linear and
continuous map $U:S(X)\rightarrow K$ given via $U(f)=\int fd\mu $ can be
extended by uniform continuity to $V:\overline{S(X)}=TM(X)\rightarrow K$.
For any $f\in TM(X)$, we write $V(f)\overset{def}{=}\int fd\mu =$ \textit{%
the integral of }$f$\textit{\ with respect to }$\mu $. Hence, for any $f\in
TM(X)$, one has $\int fd\mu =\underset{n}{\lim }\int f_{n}d\mu $, where $%
(f_{n})_{n}\subset S(X)$ is such that $f_{n}\overset{\text{u}}{\underset{n}{%
\rightarrow }}f$ (the result does not depend upon the sequence $(f_{n})_{n}$
used). So, our integral is uniform and sesquilinear ($\int (\alpha f+\beta
g)d\mu =\alpha \int fd\mu +\beta \int gd\mu $ and $\int fd(\alpha \mu +\beta
\nu )=\overline{\alpha }\int fd\mu +\overline{\beta }\int fd\nu $ for $%
\alpha $ and $\beta $ in $K$, $f$ and $g$ in $TM(X)$ and $\mu $, $\nu $ in $%
cabv(X)$). Notice that, for any $f\in TM(X)$ and any $\mu \in cabv\left(
X\right) $ one has $\left\vert \int fd\mu \right\vert \leq \left\Vert \mu
\right\Vert \left\Vert f\right\Vert $ and $\int fd\mu $ can be computed for
any $f\in C(X)\subset TM(X)$. \textit{From now on, we shall discuss about }$%
\int fd\mu $\textit{\ only for }$f\in C(X)$\textit{.}

\textit{\smallskip }

We feel obliged to insist upon some computational aspects derived from the
fact that the complex Hilbert spaces (with their sesquilinear scalar
products) make life a bit more complicated.

\textit{Firstly}, in the particular case when $X=K$, due to the fact that
the scalar product in $K$ is given by the formula $(\alpha \mid \beta
)=\alpha \overline{\beta }$, we have for any $f\in C(K)$ and any $\mu \in
cabv\left( K\right) $: $\int fd\mu $ (with the present definition) $=$ $\int
fd\overline{\mu }$ (with the standard definition, where $\overline{\mu }\in
cabv\left( K\right) $ acts via $\overline{\mu }(A)=\overline{\mu (A)}$ for
any $A\in \mathcal{B}$). Extending these considerations and considering an
orthonormal basis $\left( e_{i}\right) _{i\in I}$ of $X$, one can identify
any $f\in C(X)$ via $f\equiv \left( f_{i}\right) _{i\in I}\subset C(K)$ and
any $\mu \in cabv\left( X\right) $ via $\mu \equiv \left( \mu _{i}\right)
_{i\in I}\subset cabv\left( K\right) $, with the following explanations: a)
For any $f\in C(X)$ and any $t\in T$, $f\left( t\right) =\underset{i}{\emph{S%
}}f_{i}(t)e_{i}$. b) For any $\mu \in cabv\left( X\right) $ and any $A\in 
\mathcal{B}$, $\mu (A)=\underset{i}{\emph{S}}\mu _{i}(A)e_{i}$ (summable
families). Then one can prove that $\int fd\mu =\underset{i}{\emph{S}}\int
f_{i}d\mu _{i}$ (the integral being computed with the present definition).

\textit{Secondly}, we recall the Riesz-Fr\'{e}chet representation theorem
asserting the existence of the antilinear and isometric bijection $%
F:X\rightarrow X^{\prime }$, given via $F(y)=T_{y}$, where $T_{y}(x)=(x\mid
y)$ for any $y\in X$ and any $x\in X$. On the basis of this representation
theorem, we interpret a (now) classical result of N. Dinculeanu ([5]) and
obtain an antilinear and isometric isomorphism $H:cabv(X)\rightarrow
C(X)^{^{\prime }}$ given via $H(\mu )=V_{\mu }$, where $V_{\mu }(f)=\int
fd\mu $ for any $\mu \in cabv\left( X\right) $ and any $f\in C(X)$.

We add some more computational facts. Let $Y\subset X$ be a closed linear
subspace and $\pi _{Y}:X\rightarrow X$ the orthogonal projection defined by $%
Y$ ($\pi _{Y}(y)=y$ for any $y\in Y$). Let $f\in C(X)$ and $\mu \in
cabv\left( X\right) $. Assume that either $f(T)\subset Y$ or $\mu (B)\subset
Y$. Then $\int fd\mu =\int (\pi _{Y}\circ f)d(\pi _{Y}\circ \mu )$. Other
result for $f\in C(K)$, $\mu \in cabv\left( K\right) $, $x,y\in X$: one has $%
\int (fx)d(\mu y)=(\int fd\mu )\cdot (x\mid y)$ (in case $x=y$ and $%
\left\Vert x\right\Vert =1$, one has $\int (fx)d(\mu x)=(\int fd\mu )$).
Finally, we consider, for any $t\in T$, the Dirac measure concentrated at $t$%
, namely $\delta _{t}:\mathcal{B}\rightarrow K$, $\delta _{t}(A)=\varphi
_{A}(t)$. Then, for any $x\in X$ and any $t\in T$, $\delta _{t}x\in cabv(X)$
and $\left\Vert \delta _{t}x\right\Vert =\left\Vert x\right\Vert $. For any $%
f\in C(X)$, one has $\int fd(\delta _{t}x)=(f(t)\mid x)$.

For general measure theory see [9] and [17]. For vector measure and vector
integration see [5], [4] and [3].

\bigskip

\textbf{3. Results}

\bigskip

\textbf{The space }$L(X)$

\bigskip

\textit{From now on, }$(T,d)$\textit{\ will be a compact metric space such
that }$T$\textit{\ has at least two elements and }$X$\textit{\ will be a non
null Hilbert space with scalar product }$(.\mid .)$\textit{\ and
corresponding norm }$\left\Vert .\right\Vert $\textit{.}

Recall that a function $f:T\rightarrow X$ is a \textit{Lipschitz function}
if there exists a number $M\in (0,\infty )$ such that $\left\Vert
f(x)-f(y)\right\Vert \leq Md(x,y)$ for any $x$ and $y$ in $T$. For such $f$,
we define the \textit{Lipschitz constant of} $f$, denoted by $\left\Vert
f\right\Vert _{L}$, by $\left\Vert f\right\Vert _{L}=\sup \{\frac{\left\Vert
f(x)-f(y)\right\Vert }{d(x,y)}\mid x,y\in T,x\neq y\}$. Notice that $%
\left\Vert f\right\Vert _{L}=\min \{M\mid M\geq 0$ such that $\left\Vert
f(x)-f(y)\right\Vert \leq Md(x,y)$ for any $x,y\in T\}$.

The space 
\begin{equation*}
L(X)=\{f:T\rightarrow X\mid f\text{ is a Lipschitz function}\}
\end{equation*}%
is seminormed with the seminorm $f\mapsto \left\Vert f\right\Vert _{L}$ (we
have $\left\Vert f\right\Vert _{L}=0$ if and only if $f$ is constant). The
space $L(X)$ is normed with the norm $f\mapsto \left\Vert f\right\Vert _{BL}$
given via $\left\Vert f\right\Vert _{BL}=\left\Vert f\right\Vert +\left\Vert
f\right\Vert _{L}$. The unit ball of $L(X)$ is $BL_{1}(X)=\{f\in L(X)\mid
\left\Vert f\right\Vert _{BL}\leq 1\}$.

Of course $L(X)\subset C(X)$. In case $X=K$, $L(K)$ is dense in $C(K)$. This
assertion remains valid for $X=K^{n}$, namely we have

\bigskip

\textbf{Theorem 1}. \textit{For any} $n\in \mathbb{N}$\textit{, }$L(K^{n})$ 
\textit{is dense in }$C(K^{n})$\textit{. More precisely, there exists a
sequence }$(f^{m})_{m}\subset L(K^{n})$\textit{\ such that the set }$%
\{f^{m}\mid m\in \mathbb{N}\}$\textit{\ is dense in }$C(K^{n})$\textit{\
(which is separable).}

\textbf{Sketch of the proof}. Let $(g^{m})_{m}\subset L(K)$\textit{\ }be a
sequence such that\textit{\ }$A=\{g^{m}\mid m\in \mathbb{N}\}$ is dense in $%
C(K)$. Then the countable set $A^{n}$ is dense in $C(K^{n})$ and this proves
all, in view of the following two facts:

a) We have $A\subset L(K^{n})$ (because, if $f=(f_{1},f_{2},...,f_{n})\in
A^{n}$ and $x,y\in T$, one has

\begin{equation*}
\left\Vert f(x)-f(y)\right\Vert \leq \overset{n}{\underset{i=1}{\sum }}%
\left\vert f_{i}(x)-f_{i}(y)\right\vert \leq (\overset{n}{\underset{i=1}{%
\sum }}\left\Vert f_{i}\right\Vert _{L})d(x,y)\text{.}
\end{equation*}

b) For any $f=(f_{1},f_{2},...,f_{n})\in C(K^{n})$, there exist $n$
sequences $(g_{i}^{m})_{m}\subset A$, $i\in \{1,2,...,n\}$, such that $%
g_{i}^{m}\underset{m}{\overset{u}{\rightarrow }}f_{i}$, for all $i$ and $%
g^{m}\underset{m}{\overset{u}{\rightarrow }}f$, where $%
g^{m}=(g_{1}^{m},g_{2}^{m},...,g_{n}^{m})$. $\square $

\bigskip

\textbf{The Monge-Kantorovich norm}

\bigskip

In this subparagraph, we introduce the Monge-Kantorovich norm.

\bigskip

The next result is probably well-known, but we think a careful proof of it
is desirable. Besides, some technical parts of the proof will be used later.

\bigskip

\textbf{Theorem 2 (Lipschitz Urysohn-Type Lemma)}. \textit{Let }$\emptyset
\neq H\subset T$\textit{, }$T\neq D\subset T$\textit{\ such that }$H$\textit{%
\ is compact, }$D$\textit{\ is open and }$H\subset D$\textit{. Then there
exist a number }$M\in (0,1)$\textit{\ and a function }$f\in BL_{1}(\mathbb{R}%
)$\textit{\ with }$\left\Vert f\right\Vert _{BL}=1$\textit{\ such that }$%
0\leq f(t)\leq M$\textit{\ for any }$t\in T$,\textit{\ }$f(t)=M$\textit{\
for any }$t\in H$\textit{\ and }$f(t)=0$\textit{\ for any }$t\in C_{D}$%
\textit{.}

\textbf{Proof}. a) One has $C_{D}\neq \emptyset $ and we shall prove that $%
\delta (H,C_{D})>0$. Indeed, accepting that $\delta (H,C_{D})=0$, we find
the sequences $(x_{n})_{n}\subset H$ and $(y_{n})_{n}\subset C_{D}$ such
that $\underset{n\rightarrow \infty }{\lim }d(x_{n},y_{n})=0$. Due to the
compactness, we find (taking convergent subsequences) $x\in H$ and $y\in
C_{D}$ such that $d(x,y)=0$, i.e. $x=y\in H\cap C_{D}$, impossible.

b) For any $t\in T$, one has%
\begin{equation}
d_{1}(t,C_{D})+d_{1}(t,H)\geq \delta (C_{D},H)>0  \tag{1}
\end{equation}%
because, for any sequences $(a_{n})_{n}\subset C_{D}$ and $%
(b_{n})_{n}\subset H$, one has 
\begin{equation*}
d(t,a_{n})+d(t,b_{n})\geq d(a_{n},b_{n})\geq \delta (C_{D},H)\text{.}
\end{equation*}

So, one can define $g:T\rightarrow \mathbb{R}_{+}$, via%
\begin{equation*}
g(t)=\frac{d_{1}(t,C_{D})}{d_{1}(t,C_{D})+d_{1}(t,H)}
\end{equation*}%
and $0\leq g(t)\leq 1$ for any $t\in T$, $g(t)=1$ for any $t\in H$ and $%
g(t)=0$ for any $t\in C_{D}$.

For $x$ and $y$ arbitrarily taken in $T$, using $(1)$, we have:%
\begin{equation*}
\left\vert g(x)-g(y)\right\vert =\frac{\left\vert
d_{1}(x,C_{D})d_{1}(y,H)-d_{1}(x,H)d_{1}(y,C_{D})\right\vert }{%
(d_{1}(x,C_{D})+d_{1}(x,H))(d_{1}(y,C_{D})+d_{1}(y,H))}\leq
\end{equation*}%
\begin{equation*}
\leq \frac{\left\vert
d_{1}(x,C_{D})d_{1}(y,H)-d_{1}(y,C_{D})d_{1}(y,H)+d_{1}(y,C_{D})d_{1}(y,H)-d_{1}(x,H)d_{1}(y,C_{D})\right\vert 
}{(\delta (C_{D},H))^{2}}\leq
\end{equation*}%
\begin{equation*}
\leq \frac{d_{1}(y,H)\left\vert d_{1}(x,C_{D})-d_{1}(y,C_{D})\right\vert
+d_{1}(y,C_{D})\left\vert d_{1}(y,H)-d_{1}(x,H)\right\vert }{(\delta
(C_{D},H))^{2}}\leq
\end{equation*}%
\begin{equation*}
\leq \frac{(d_{1}(y,H)+d_{1}(y,C_{D}))d(x,y)}{(\delta (C_{D},H))^{2}}\leq 
\frac{2diam(T)}{(\delta (C_{D},H))^{2}}d(x,y)
\end{equation*}%
The last inequality is true because $H$ and $C_{D}$ being compact, we have $%
d_{1}(y,H)=d(y,h)$ for some $h\in H$ and $d_{1}(y,C_{D})=d(y,p)$ for some $%
p\in C_{D}$ a.s.o.

It follows that $g$ is a Lipschitz function and for any $x$ and $y$ in $T$
one has $\left\vert g(x)-g(y)\right\vert \leq Bd(x,y)$, where ($T$ has at
least two points) $B=\frac{2diam(T)}{(\delta (C_{D},H))^{2}}>0$.

Hence $\left\Vert g\right\Vert =1$, $\left\Vert g\right\Vert _{L}\leq B$, $g$
is not constant and consequently $1+B\geq 1+\left\Vert g\right\Vert
_{L}=\left\Vert g\right\Vert _{BL}>1$.

Finally, we define 
\begin{equation*}
f=\frac{1}{\left\Vert g\right\Vert _{BL}}g\text{ and }M=\frac{1}{\left\Vert
g\right\Vert _{BL}}\text{. }\square
\end{equation*}

\bigskip

\textbf{Lemma 3}. \textit{Let }$\mu _{1}:\mathcal{B}\rightarrow \mathbb{R}%
_{+}$\textit{\ and }$\mu _{2}:\mathcal{B}\rightarrow \mathbb{R}_{+}$\textit{%
\ be two finite }$\sigma $-\textit{additive measures. We have the
equivalence:} $\mu _{1}=\mu _{2}\Longleftrightarrow $ $\int fd\mu _{1}=\int
fd\mu _{2}$ \textit{for any positive function} $f\in BL_{1}(\mathbb{R})$.

\textbf{Proof}. One must prove the implication $"\Leftarrow "$.

Because $T$ is a metric space, the measures $\mu _{1}$ and $\mu _{2}$ are
regular: for any $A\in \mathcal{B}$ one has $\mu _{1}(A)=\sup \mu _{1}(H)$
and $\mu _{2}(A)=\sup \mu _{2}(H)$, the suprema being computed for all
compact subsets $H\subset A$. So, it will suffice to prove that $\mu
_{1}(H)=\mu _{2}(H)$ for any compact subset $H\subset T$. Take such a
compact $H$.

\textit{First Possibility}: the only open subset $D$ of $T$ such that $%
H\subset D$ is $D=T$. Again the regularity of $\mu _{1}$ and $\mu _{2}$ says
that $\mu _{i}(H)=\inf \{\mu _{i}(D)\mid H\subset D\subset T$, $D$ open$\}$,
hence $\mu _{i}(H)=\mu _{i}(T)$, for any $i\in \{1,2\}$.

Taking $f:T\rightarrow \mathbb{R}$, $f(t)=1$ for any $t\in T$, one has $f\in
BL_{1}(\mathbb{R})$ and $\int fd\mu _{1}=\int fd\mu _{2}$ (according to the
hypothesis), hence $\mu _{1}(T)=\mu _{2}(T)$, i.e. $\mu _{1}(H)=\mu _{2}(H)$.

\textit{Second Possibility}: there exist an open set $\Delta $ such that $%
H\subset \Delta \subset T$, $\Delta \neq T$. We shall prove that $\mu
_{1}(H)\leq \mu _{2}(H)$ (and, in the same way, $\mu _{2}(H)\leq \mu _{1}(H)$%
), hence $\mu _{1}(H)=\mu _{2}(H)$.

Take arbitrarily $\varepsilon >0$. The regularity of $\mu _{2}$ yields an
open set $D_{\varepsilon }$ such that $H\subset D_{\varepsilon }$ and $\mu
_{2}(D_{\varepsilon })\leq \mu _{2}(H)+\varepsilon $. Write $D=\Delta \cap
D_{\varepsilon }$, hence $D$ is open, $C_{D}\neq \emptyset $ and $\mu
_{2}(D)\leq \mu _{2}(D_{\varepsilon })\leq \mu _{2}(H)+\varepsilon $. For
the couple $(H,D)$, we construct the function $f$ from Theorem 2. According
to the hypothesis we have $\int fd\mu _{1}=\int fd\mu _{2}$ and this implies%
\begin{equation*}
M\mu _{1}(H)\leq \int fd\mu _{1}=\int fd\mu _{2}\leq M\mu _{2}(D)\leq M(\mu
_{2}(H)+\varepsilon )\text{,}
\end{equation*}%
so $\mu _{1}(H)\leq \mu _{2}(H)+\varepsilon $. Because $\varepsilon $ is
arbitrary, it follows that $\mu _{1}(H)\leq \mu _{2}(H)$. $\square $

\bigskip

\textbf{Theorem 4}. \textit{For any }$\mu \in cabv(X)$, \textit{one has the
equivalence:} $\mu =0$\textit{\ }$\Leftrightarrow $\textit{\ }$\int fd\mu =0$
\textit{for any} $f\in L(X)$.

\textbf{Proof}. We must prove the implication $"\Leftarrow "$.

\textit{Case }$X=\mathbb{R}$. Write $\mu =\mu _{1}-\mu _{2}$, where $\mu
_{1},\mu _{2}\geq 0$ and $\mu _{1},\mu _{2}\in cabv(\mathbb{R})$ (Jordan
decomposition). Hence, for any $f\in L(\mathbb{R})$, $0=\int fd\mu =\int
fd\mu _{1}-\int fd\mu _{2}$. Using Lemma 3, we get $\mu _{1}=\mu _{2}$,
hence $\mu =0$.

\textit{Case} $X=\mathbb{C}$. Write $\mu =\mu _{1}+i\mu _{2}$, where $\mu
_{1}$ and $\mu _{2}$ are in $cabv(\mathbb{R})$. Hence, for any $f\in L(%
\mathbb{R})$, one has $0=\int fd\mu =\int fd\mu _{1}+i\int fd\mu
_{2}\Leftrightarrow \int fd\mu _{1}=\int fd\mu _{2}=0$. Using the case $X=%
\mathbb{R}$, we get $\mu _{1}=\mu _{2}=0$, hence $\mu =0$.

\textit{General case}. Let $(e_{i})_{i\in I}$ be an orthonormal basis for $X$%
. Write, for any $x\in X$: $x=\underset{i\in I}{\emph{S}}a_{i}e_{i}$, hence $%
\left\Vert x\right\Vert ^{2}=\underset{i\in I}{\emph{S}}\left\vert
a_{i}\right\vert ^{2}$ (summable families, $a_{i}=(x\mid e_{i})$ Fourier
coefficients). For any $i\in I$, define $H_{i}:X\rightarrow K$ via $%
H_{i}(x)=a_{i}$. Then $H_{i}\in X^{\prime }$ and $\left\Vert
H_{i}\right\Vert _{0}=1$.

For any $\mu \in cabv(X)$ and any $A\in \mathcal{B}$, one has $\mu (A)=%
\underset{i\in I}{\emph{S}}\mu _{i}(A)e_{i}$, where $\mu _{i}=H_{i}\circ \mu
\in cabv(K)$ are uniquely determined. The fact that all $\mu _{i}$ are $%
\sigma $-additive is obvious, due to the continuity of $H_{i}$. It follows
that all $\mu _{i}$ are of bounded variation (actually $\left\Vert \mu
_{i}\right\Vert \leq \left\Vert \mu \right\Vert $ for any $i\in I$, but we
will not use this fact).

Now, take $\mu \in cabv(X)$ such that $\int fd\mu =0$ for any $f\in L(X)$.
To prove that $\mu =0$ means to prove that all $\mu _{i}=0$. Fix $i\in I$
arbitrarily. For any $\varphi \in L(K)$, one has $\varphi e_{i}\in L(X)$ and 
$\varphi e_{i}(T)\subset Sp(e_{i})$, $(\mu _{i}e_{i})(\mathcal{B})\subset
Sp(e_{i})$. Let $Y_{i}=Sp(e_{i})$. Obviously $\varphi e_{i}=\pi
_{Y_{i}}\circ (\varphi e_{i})$ and $\mu _{i}e_{i}=\pi _{Y_{i}}\circ \mu $.
According to the hypothesis: 
\begin{equation*}
0=\int (\varphi e_{i})d\mu =\int (\pi _{Y_{i}}\circ (\varphi e_{i}))d\mu
=\int (\pi _{Y_{i}}\circ (\varphi e_{i}))d(\pi _{Y_{i}}\circ \mu )=
\end{equation*}%
\begin{equation*}
=\int (\varphi e_{i})d(\mu _{i}e_{i})=(\int \varphi d\mu _{i})\left\Vert
e_{i}\right\Vert ^{2}=\int \varphi d\mu _{i}\text{.}
\end{equation*}

Because $\varphi \in L(K)$ is arbitrary, we get $\mu _{i}=0$. $\square $

\bigskip

\textbf{Theorem 5}. \textit{For any }$\mu \in cabv(X)$\textit{, define}%
\begin{equation*}
\left\Vert \mu \right\Vert _{MK}=\sup \{\left\vert \int fd\mu \right\vert
\mid f\in BL_{1}(X)\}\text{.}
\end{equation*}

\textit{Then, the function} $\mu \mapsto \left\Vert \mu \right\Vert _{MK}$ 
\textit{is a norm on} $cabv(X)$ \textit{and one has}%
\begin{equation*}
\left\Vert \mu \right\Vert _{MK}\leq \left\Vert \mu \right\Vert
\end{equation*}%
\textit{for any} $\mu \in cabv(X)$.

\textbf{Proof}. Let $\mu \in cabv(X)$. Then, in view of the antilinear
identification $C(X)^{^{\prime }}\equiv cabv(X)$, one has%
\begin{equation*}
\left\Vert \mu \right\Vert =\sup \{\left\vert \int fd\mu \right\vert \mid
f\in B^{1}(X)\}\text{,}
\end{equation*}%
where $B^{1}(X)=\{f\in C(X)\mid \left\Vert f\right\Vert \leq 1\}$ and $%
\left\Vert f\right\Vert \leq \left\Vert f\right\Vert _{BL}$ implies $%
BL_{1}(X)\subset B^{1}(X)$, hence $\left\Vert \mu \right\Vert _{MK}\leq
\left\Vert \mu \right\Vert $.

It is obvious that $\mu \mapsto \left\Vert \mu \right\Vert _{MK}$ is a
seminorm. To finish the proof, one must show the implication $\left\Vert \mu
\right\Vert _{MK}=0\Rightarrow \mu =0$.

But $\left\Vert \mu \right\Vert _{MK}=0$ means $\int fd\mu =0$ for any $f\in
L(X)$ and this implies $\mu =0$, according to Theorem 4. $\square $

\bigskip

\textbf{Definition 6}. The norm $\left\Vert .\right\Vert _{MK}$\ is called
the \textit{Monge-Kantorovich norm}.

\bigskip

Notice that, according to the definition, one has, for any $\mu \in cabv(X)$
and any $f\in L(X)$:%
\begin{equation}
\left\vert \int fd\mu \right\vert \leq \left\Vert \mu \right\Vert
_{MK}\left\Vert f\right\Vert _{BL}\text{.}  \tag{2}
\end{equation}

For any $t\in T$ and any $x\in X$, $\left\Vert x\right\Vert =1$, one has%
\begin{equation}
\left\Vert \delta _{t}x\right\Vert _{MK}=1\text{.}  \tag{3}
\end{equation}

Indeed, write $\delta _{t}x=\mu $ and take $f\in BL_{1}(X)$. One has $\int
fd\mu =(f(t)\mid x)$, hence $\left\vert \int fd\mu \right\vert \leq
\left\Vert f(t)\right\Vert \left\Vert x\right\Vert =\left\Vert
f(t)\right\Vert \leq 1$ and $\left\Vert \mu \right\Vert _{MK}\leq 1$.
Conversely, taking $f(t)=x$ for any $t\in T$, one has $\left\Vert
f\right\Vert _{BL}=\left\Vert f\right\Vert =1$, hence $\left\vert \int fd\mu
\right\vert =(x\mid x)=1$ and $\left\Vert \mu \right\Vert _{MK}\geq 1$ a.s.o.

The topology generated by $\left\Vert .\right\Vert _{MK}$ on $cabv(X)$ will
be denoted by $\mathcal{T}(MK,X)$ (\textit{the Monge-Kantorovich topology}).
For any $a>0$, the topology induced by $\mathcal{T}(MK,X)$ on $B_{a}(X)$
will be denoted by $\mathcal{T}(MK,X,a)$. For a sequence $(\mu
_{n})_{n}\subset cabv(X)$ and for $\mu \in cabv(X)$, we shall write $\mu _{n}%
\underset{n}{\overset{\text{MK}}{\rightarrow }}\mu $ to denote the fact that 
$(\mu _{n})_{n}$ converges to $\mu $ in the Monge-Kantorovich topology.

In the sequel, we shall make some considerations concerning the comparison
between the variational topology $\mathcal{T}(var,X)$ and the
Monge-Kantorovich topology $\mathcal{T}(MK,X)$.

Due to the inequality $\left\Vert \mu \right\Vert _{MK}\leq \left\Vert \mu
\right\Vert $, we have $\mathcal{T}(MK,X)\subset \mathcal{T}(var,X)$. Of
course, if $T$ is finite, one has $\mathcal{T}(MK,X)=\mathcal{T}(var,X)$. As
concerns the case when $T$ is infinite, we remark first that $T$ is infinite
if and only if $T$ has at least an accumulation point. Here comes

\bigskip

\textbf{Theorem 7}. \textit{Assume }$T$\textit{\ is infinite. Then the
inclusion }$\mathcal{T}(MK,X)\subset \mathcal{T}(var,X)$\textit{\ is strict.
Also in this case, the normed space }$(cabv(X),\left\Vert .\right\Vert
_{MK}) $ \textit{is not Banach}.

\textbf{Proof}. a) First we shall prove that for any $a$ and $b$ in $T$, $%
a\neq b$ and any $x\in X$, $\left\Vert x\right\Vert =1$, one has 
\begin{equation}
\left\Vert \delta _{a}x-\delta _{b}x\right\Vert =2  \tag{4}
\end{equation}%
(in case $X=K$, $x=1$, one has $\left\Vert \delta _{a}-\delta
_{b}\right\Vert =2$).

To prove $(4)$, write $\delta _{a}x-\delta _{b}x=\mu $ and take a partition $%
(A_{i})_{i\in \{1,2,...,m\}}$ of $T$. One has either $\overset{m}{\underset{%
i=1}{\sum }}\left\Vert \mu (A_{i})\right\Vert =0$ (in case there exists $i$
such that $a\in A_{i}$ and $b\in A_{i}$) or $\overset{m}{\underset{i=1}{\sum 
}}\left\Vert \mu (A_{i})\right\Vert =2$ (in case there exist $i\neq j$ such
that $a\in A_{i}$ and $b\in A_{j}$). The second alternative is always
possible, taking $A_{1}=B(a,r)$, $A_{2}=B(b,r)$ with $A_{1}\cap
A_{2}=\emptyset $ and the other $A_{i}$ arbitrarily. Hence, by passing to
supremum, one gets $\left\Vert \mu \right\Vert =2$.

b) Again, for $a$ and $b$ in $T$, $a\neq b$ and any $x\in X$, $\left\Vert
x\right\Vert =1$, we shall prove that 
\begin{equation}
\left\Vert \delta _{a}x-\delta _{b}x\right\Vert _{MK}\leq d(a,b)\text{.} 
\tag{5}
\end{equation}

Indeed, writing again $\mu =\delta _{a}x-\delta _{b}x$, we have, for any $%
f\in L(X)$: $\int fd\mu =(f(a)-f(b)\mid x)$ (in case $X=K$, $x=1$: $\int
fd\mu =f(a)-f(b)$). Hence, if $f\in BL_{1}(X)$, one has%
\begin{equation*}
\left\vert \int fd\mu \right\vert \leq \left\Vert f(a)-f(b)\right\Vert
\left\Vert x\right\Vert =\left\Vert f(a)-f(b)\right\Vert \leq d(a,b)
\end{equation*}%
and passing to supremum, we get $(5)$.

At the end of the paper, we shall discuss supplementarily formula $(5)$.

c) Because $T$ is infinite, we take an accumulation point $t_{0}\in T$ and a
sequence $(t_{n})_{n}\subset T$ such that $t_{n}\underset{n}{\rightarrow }%
t_{0}$ and $t_{n}\neq t_{0}$ for any $n\geq 1$. According to $(5)$, it
follows that, for any $x\in X$ with $\left\Vert x\right\Vert =1$, one has $%
\delta _{t_{n}}x\underset{n}{\overset{\text{MK}}{\rightarrow }}\delta
_{t_{0}}x$, whereas, according to $(4)$, the assertion $\delta _{t_{n}}x%
\underset{n}{\overset{\text{var}}{\rightarrow }}\delta _{t_{0}}x$ is false.
Hence the inclusion $\mathcal{T}(MK,X)\subset \mathcal{T}(var,X)$ must be
strict.

The fact that $(cabv(X),\left\Vert .\right\Vert _{MK})$ is not Banach
follows from the inequality $\left\Vert .\right\Vert _{MK}\leq \left\Vert
.\right\Vert $ and from the fact that the norms $\left\Vert .\right\Vert
_{MK}$ and $\left\Vert .\right\Vert $ are not equivalent. $\square $

\bigskip

Let us present some Supplementary Remarks

\bigskip

\textbf{Remarks}

a) If $T$ is infinite, one can find a sequence $(\mu _{n})_{n}\subset
cabv(X) $ such that $\left\Vert \mu _{n}\right\Vert _{MK}=1$ and $\left\Vert
\mu _{n}\right\Vert >n$ for any $n$.

b) Generally speaking one has for any sequence $(\mu _{n})_{n}\subset
cabv(X) $ the implication $\mu _{n}\underset{n}{\overset{\text{var}}{%
\rightarrow }}\mu \Rightarrow \mu _{n}\underset{n}{\overset{\text{MK}}{%
\rightarrow }}\mu $. The converse implication is not true for infinite $T$.

\bigskip

\textbf{Example}

Take $T=[0,1]$, $X=\mathbb{R}$ and consider the true fact that $\delta _{%
\frac{1}{n}}\underset{n}{\overset{\text{MK}}{\rightarrow }}\delta _{0}$.
Because, for any $n\in \mathbb{N}$, one has $\delta _{\frac{1}{n}}((0,1])=1$
and $\delta _{0}((0,1])=0$, it follows that $\delta _{\frac{1}{n}}\underset{n%
}{\overset{\text{var}}{\rightarrow }}\delta _{0}$ is a false assertion
because, generally speaking $\mu _{n}\underset{n}{\overset{\text{var}}{%
\rightarrow }}\mu $ implies $\mu _{n}(A)\underset{n}{\rightarrow }\mu (A)$
for any $A\in \mathcal{B}$.

So, \textit{convergence in the Monge-Kantorovich norm does not imply
pointwise convergence}. We shall see later that convergence in the
Monge-Kantorovich topology means weak$^{\ast }$ convergence for bounded
sequences and this explains everything.

\bigskip

\textbf{The weak}$^{\ast }$\textbf{\ topology on }$cabv(X)$

\bigskip

Let us introduce a new topology on $cabv(X)$. This topology is defined on
the basis of the fact that $cabv(X)$ is identified with the dual of $C(X)$.

\bigskip

\textbf{Definition 8}. \textit{The weak}$^{\ast }$\textit{\ topology on} $%
cabv(X)$ is the (separated) locally convex topology on $cabv(X)$ generated
by the family of seminorms $(p_{f})_{f\in C(X)}$, where, for any $f\in C(X)$%
, $p_{f}:cabv(X)\rightarrow \mathbb{R}_{+}$ is given via%
\begin{equation*}
p_{f}(\mu )=\left\vert \int fd\mu \right\vert \text{.}
\end{equation*}

\bigskip

The weak$^{\ast }$ topology will be denoted by $\mathcal{T(}w^{\ast },X%
\mathcal{)}$ and, for any $a>0$, its restriction to $B_{a}(X)$, will be
denoted by $\mathcal{T(}w^{\ast },X,a\mathcal{)}$.

For any $\mu \in cabv(X)$, a neighborhood basis for $\mu $ is formed with
all sets of the form 
\begin{equation*}
V(\mu ;g_{1},g_{2},...,g_{m};\varepsilon )\overset{def}{=}\{\upsilon \in
cabv(X)\mid \left\vert \int g_{i}d(\mu -\upsilon )\right\vert <\varepsilon 
\text{, }i\in \{1,2,...,m\}\}
\end{equation*}
(one takes into consideration all possible $\varepsilon >0$, all $m\in 
\mathbb{N}$ and all $g_{i}\in C(X)$).

For a sequence $(\mu _{n})_{n}\subset cabv(X)$ and for $\mu \in cabv(X)$, we
shall write $\mu _{m}\overset{\text{w}^{\ast }}{\underset{m}{\rightarrow }}%
\mu $ to denote the fact that $(\mu _{n})_{n}$ converges to $\mu $ in $%
\mathcal{T(}w^{\ast },X\mathcal{)}$. This means that $\underset{m}{\lim }%
\int fd\mu _{m}=\int fd\mu $ for any $f\in C(X)$.

Notice that Alaoglu's theorem implies that, for any $a>0$, the set $B_{a}(X)$
is weak$^{\ast }$ compact (i.e. compact in $\mathcal{T(}w^{\ast },X\mathcal{)%
}$).

In the sequel, we shall fix $n\in \mathbb{N}$ and we shall work for $X=K^{n}$%
. \ We have seen (Theorem 1) that one can find a sequence $%
(f_{m})_{m}\subset L(K^{n})$ such that $\{f_{m}\mid m\in \mathbb{N}\}$ is
dense in $C(K^{n})$. This fact has the following two important consequences.

\bigskip

\textbf{Theorem 9 (Metrisability of }$B_{a}(K^{n})$\textbf{\ under }$%
\mathcal{T}(w^{\ast },K^{n})$\textbf{)}. \textit{For any }$a>0$\textit{, the
topology }$\mathcal{T}(w^{\ast },K^{n},a)$ \textit{is metrisable. The set }$%
B_{a}(K^{n})$ \textit{is compact as a subset of the topological space} $%
(cabv(K^{n}),\mathcal{T}(w^{\ast },K^{n})$\textit{). Consequently,} $%
B_{a}(K^{n})$ \textit{considered as a metric space (with any metric
generating} $\mathcal{T}(w^{\ast },K^{n},a)$\textit{) is complete.}

\bigskip

The metrisability of $\mathcal{T}(w^{\ast },K^{n},a)$ follows from the
separability of $C(K^{n})$, viewing $cabv(K^{n})$ as the dual of $C(K^{n})$
(see [8], V, 5.1, page 426).

\bigskip

\textbf{Theorem 10}. \textit{Let} $(f_{m})_{m}\subset L(K^{n})$ \textit{be
the aforementioned dense sequence in }$C(K^{n})$\textit{. Then, for any }$%
a>0 $\textit{, any sequence }$(\mu _{p})_{p}\subset B_{a}(K^{n})$\textit{\
and any }$\mu \in B_{a}(K^{n})$\textit{\ one has the equivalence: }$\mu _{p}%
\overset{w^{\ast }}{\underset{p}{\rightarrow }}\mu $\textit{\ }$%
\Leftrightarrow $\textit{\ }$\int f_{m}d\mu _{p}\underset{p}{\rightarrow }%
\int f_{m}d\mu $\textit{\ for any} $m\in \mathbb{N}$.

\textbf{Proof}. Only the implication $"\Leftarrow "$ must be proved. Let us
consider $V=V(\mu ;g_{1},g_{2},...,g_{p};\varepsilon )\cap B_{a}(X)$ a basic
neighborhood of $\mu $ in $\mathcal{T(}w^{\ast },X,a\mathcal{)}$. For any $%
m\in \{1,2,...,p\}$, choose $f_{i_{m}}$ such that $\left\Vert
f_{i_{m}}-g_{m}\right\Vert _{\infty }<\frac{\varepsilon }{4a}$. Take $\delta
=\min \{\frac{\varepsilon }{3},\frac{\varepsilon }{4a}\}$, construct $%
W=V(\mu ;f_{i_{1}},f_{i_{2}},...,f_{i_{p}};\delta )\cap B_{a}(X)$ and notice
that $W\subset V$. Indeed, if $\upsilon \in W$, one has, for any $m\in
\{1,2,...,p\}$:%
\begin{equation*}
\left\vert \int g_{m}d(\mu -\upsilon )\right\vert \leq \left\vert \int
(g_{m}-f_{i_{m}})d(\mu -\upsilon )\right\vert +\left\vert \int
f_{i_{m}}d(\mu -\upsilon )\right\vert \leq
\end{equation*}%
\begin{equation*}
\leq \left\Vert g_{m}-f_{i_{m}}\right\Vert \left\vert \mu -\upsilon
\right\vert (T)+\delta \leq \frac{\varepsilon }{4a}\left\Vert \mu -\upsilon
\right\Vert +\frac{\varepsilon }{3}\leq
\end{equation*}%
\begin{equation*}
\leq \frac{\varepsilon }{4a}(\left\Vert \mu \right\Vert +\left\Vert \upsilon
\right\Vert )+\frac{\varepsilon }{3}\leq \frac{\varepsilon }{4a}(a+a)+\frac{%
\varepsilon }{3}<\varepsilon \text{.}
\end{equation*}

Using the hypothesis, one can find $n_{V}\in \mathbb{N}$ such that $\mu
_{n}\in W\subset V$ for any $n\in \mathbb{N}$, $n\geq n_{V}$. $\square $

\bigskip

We shall need:

\bigskip

\textbf{Theorem 11} \textbf{(Arzela-Ascoli-Type Theorem)}. \textit{For any} $%
n\in \mathbb{N}$\textit{, the set} $BL_{1}(K^{n})$ \textit{is relatively
compact in }$C(K^{n})$.

\textbf{Proof}. Let $(f^{m})_{m}$ be a sequence in $BL_{1}(K^{n})$, with
each $f^{m}=(f_{1}^{m},f_{2}^{m},...,f_{n}^{m})$, $f_{i}^{m}\in C(K)$. Due
to the fact that for $x=(x_{1},x_{2},....,x_{n})\in K^{n}$ one has $%
\left\Vert x\right\Vert \geq \left\vert x_{i}\right\vert $, we get $%
\left\Vert f_{i}^{m}\right\Vert \leq \left\Vert f^{m}\right\Vert \leq
\left\Vert f^{m}\right\Vert _{BL}\leq 1$ and $\left\Vert
f_{i}^{m}\right\Vert _{L}\leq \left\Vert f^{m}\right\Vert _{L}\leq
\left\Vert f^{m}\right\Vert _{BL}\leq 1$ for any $m\in \mathbb{N}$ and $i\in
\{1,2,...,n\}$. Hence the sequence $(f_{1}^{m})_{m}$ is bounded and
equicontinuous in $C(K)$. Using the Arzela-Ascoli Theorem, we find a
subsequence $(f_{1}^{m_{1}^{p}})_{p}\subset (f_{1}^{m})_{m}$ and a function $%
f_{1}\in C(K)$ such that $f_{1}^{m_{1}^{p}}\underset{p}{\overset{u}{%
\rightarrow }}f_{1}$. Continuing, we find $(f_{2}^{m_{2}^{p}})_{p}\subset
(f_{2}^{m_{1}^{p}})_{p}$ and $f_{2}\in C(K)$ such that $f_{2}^{m_{2}^{p}}%
\underset{p}{\overset{u}{\rightarrow }}f_{2}$ and so on. Finally, we find $%
(f_{n}^{m_{n}^{p}})_{p}\subset (f_{n}^{m_{n-1}^{p}})_{p}$ and $f_{n}\in C(K)$
such that $f_{n}^{m_{n}^{p}}\underset{p}{\overset{u}{\rightarrow }}f_{n}$.

It follows that $(f^{m_{n}^{p}})_{p}\subset (f^{m})_{m}$ and $f^{m_{n}^{p}}%
\underset{p}{\overset{u}{\rightarrow }}f=(f_{1},f_{2},...,f_{n})\in C(K^{n})$%
. $\square $

\bigskip

\textbf{Remark}. It is natural to ask whether the previous result remains
valid for an arbitrary Hilbert space $X$ instead of $K^{n}$ (i.e. if $%
BL_{1}(X)$ is relatively compact in $C(X)$ also for infinite dimensional $X$%
). The answer is negative, as we shall see later.

\bigskip

We begin the investigation of the connection between the topologies $%
\mathcal{T(}w^{\ast },K^{n}\mathcal{)}$ and $\mathcal{T(}MK,K^{n}\mathcal{)}$%
.

\bigskip

\textbf{Theorem 12} \textbf{(Coincidence of weak}$^{\ast }$\textbf{%
-Convergence and Monge-Kantorovich Convergence)}. \textit{Let} $a>0$ \textit{%
and} $n\in \mathbb{N}$. \textit{For a sequence} $(\mu _{m})_{m}\subset
B_{a}(K^{n})$ \textit{and for} $\mu \in B_{a}(K^{n})$\textit{\ we have the
equivalence: }$\mu _{m}\underset{m}{\overset{\text{MK}}{\rightarrow }}\mu $%
\textit{\ }$\Leftrightarrow $ $\mu _{m}\underset{m}{\overset{\text{w}^{\ast }%
}{\rightarrow }}\mu $.

\textbf{Proof}. \textit{The implication }$"\Rightarrow "$\textit{\ }Accept
that $\mu _{m}\underset{m}{\overset{\text{MK}}{\rightarrow }}\mu $. In view
of Theorem 10, one must prove that, for any $p\in \mathbb{N}$ one has $%
\underset{m}{\lim }\int f_{p}d\mu _{m}=\int f_{p}d\mu $, where $%
(f_{p})_{p}\subset L(K^{n})$ is a dense sequence in $C(K^{n})$.

So take an arbitrary $f_{p}\neq 0$ and let $g=\frac{1}{\alpha }f_{p}$, where 
$\alpha =\left\Vert f_{p}\right\Vert _{BL}$, hence $\left\Vert g\right\Vert
_{BL}=1$. Take also $\varepsilon >0$ arbitrarily. Our hypothesis being that 
\begin{equation*}
\underset{m}{\lim }(\sup \{\left\vert \int hd(\mu _{m}-\mu )\right\vert \mid
h\in BL_{1}(K^{n})\})=0\text{,}
\end{equation*}%
one can find $m_{\varepsilon }\in \mathbb{N}$ such that for any $m\in 
\mathbb{N}$, $m\geq m_{\varepsilon }$ and any $h\in BL_{1}(K^{n})$, one has%
\begin{equation*}
\left\vert \int hd(\mu _{m}-\mu )\right\vert <\frac{\varepsilon }{\alpha }%
\text{,}
\end{equation*}%
hence%
\begin{equation*}
\left\vert \int gd(\mu _{m}-\mu )\right\vert <\frac{\varepsilon }{\alpha }%
\text{,}
\end{equation*}%
which means%
\begin{equation*}
\left\vert \int f_{p}d(\mu _{m}-\mu )\right\vert <\varepsilon \text{.}
\end{equation*}

\textit{The implication }$"\Leftarrow "$\textit{\ }Accept (reductio ad
absurdum) the existence of a sequence $(\mu _{m})_{m}\subset B_{a}(K^{n})$
and of a $\mu \in B_{a}(K^{n})$ such that $\mu _{m}\underset{m}{\overset{%
\text{w}^{\ast }}{\rightarrow }}\mu $ and such that the assertion $\mu _{m}%
\underset{m}{\overset{\text{MK}}{\rightarrow }}\mu $ is false. We shall
arrive at a contradiction.

Indeed, there exists $\varepsilon _{0}>0$ and a subsequence $(\mu
_{m_{p}})_{m_{p}}\subset (\mu _{m})_{m}$ such that $\left\Vert \mu
_{m_{p}}-\mu \right\Vert _{MK}>2\varepsilon _{0}$ for any $p$. So, for any $%
p $, one can find $f_{p}\in BL_{1}(K^{n})$ such that%
\begin{equation}
\left\vert \int f_{p}d(\mu _{m_{p}}-\mu )\right\vert >\varepsilon _{0}\text{.%
}  \tag{6}
\end{equation}

Using Theorem 11, one can find $(f_{p_{q}})_{q}\subset (f_{p})_{p}$ and $%
f\in C(K^{n})$ such that $\left\Vert f_{p_{q}}-f\right\Vert \underset{q}{%
\rightarrow }0$. Because $\mu _{m_{p}}\underset{p}{\overset{\text{w}^{\ast }}%
{\rightarrow }}\mu $, we get $p_{1}\in \mathbb{N}$ such that, for any $p\geq
p_{1}$, one has%
\begin{equation}
\left\vert \int fd(\mu _{m_{p}}-\mu )\right\vert <\frac{\varepsilon _{0}}{2}%
\text{.}  \tag{7}
\end{equation}

Let $q_{1}\in \mathbb{N}$ be such that $p_{q_{1}}>p_{1}$, and, for any $%
q\geq q_{1}$, one has%
\begin{equation}
\left\Vert f_{p_{q}}-f\right\Vert <\frac{\varepsilon _{0}}{4a}\text{.} 
\tag{8}
\end{equation}

From (6), it follows that, for any $q\geq q_{1}$, one has%
\begin{equation}
\left\vert \int f_{p_{q}}d(\mu _{m_{p_{q}}}-\mu )\right\vert >\varepsilon
_{0}\text{.}  \tag{9}
\end{equation}

At the same time, for such $q$, one has%
\begin{equation*}
\left\vert \int f_{p_{q}}d(\mu _{m_{p_{q}}}-\mu )\right\vert \leq \left\vert
\int (f_{p_{q}}-f)d(\mu _{m_{p_{q}}}-\mu )\right\vert +\left\vert \int
fd(\mu _{m_{p_{q}}}-\mu )\right\vert \leq
\end{equation*}%
\begin{equation*}
\leq \left\Vert f_{p_{q}}-f\right\Vert \left\Vert \mu _{m_{p_{q}}}-\mu
\right\Vert +\left\vert \int fd(\mu _{m_{p_{q}}}-\mu )\right\vert \leq
\end{equation*}%
\begin{equation*}
\leq \left\Vert f_{p_{q}}-f\right\Vert (\left\Vert \mu
_{m_{p_{q}}}\right\Vert +\left\Vert \mu \right\Vert )+\left\vert \int fd(\mu
_{m_{p_{q}}}-\mu )\right\vert \leq
\end{equation*}%
\begin{equation*}
\leq 2a\left\Vert f_{p_{q}}-f\right\Vert +\left\vert \int fd(\mu
_{m_{p_{q}}}-\mu )\right\vert <\frac{\varepsilon _{0}}{4a}2a+\frac{%
\varepsilon _{0}}{2}=\varepsilon _{0}\text{,}
\end{equation*}%
where we used $(7)$ and $(8)$. This contradicts $(9)$. $\square $

\bigskip

Let us interpret the last results. Take arbitrarily $a>0$ and $n\in \mathbb{N%
}$. On $B_{a}(K^{n})$ we have two metrisable topologies: $\mathcal{T(}%
MK,K^{n},a\mathcal{)}$ and $\mathcal{T(}w^{\ast },K^{n},a\mathcal{)}$ (with
Theorem 9). Theorem 12 says that the convergent sequences coincide in these
topologies, hence they are equal:

\begin{equation}
\mathcal{T(}MK,K^{n},a\mathcal{)}=\mathcal{T(}w^{\ast },K^{n},a\mathcal{)}%
\text{.}  \tag{10}
\end{equation}

Again Theorem 9 says that $B_{a}(K^{n})$ is compact for $\mathcal{T(}w^{\ast
},K^{n},a\mathcal{)}$, hence for $\mathcal{T(}MK,K^{n},a\mathcal{)}$. So $%
B_{a}(K^{n})$ is a compact (hence complete) metric space for the metric
given by $\left\Vert .\right\Vert _{MK}$.

We got (see Theorem 7 too):

\bigskip

\textbf{Theorem 13}. \textit{For any }$a>0$\textit{\ and any }$n\in \mathbb{N%
}$\textit{, the set }$B_{a}(K^{n})$\textit{, equipped with the metric
generated by the Monge-Kantorovich norm }$\left\Vert .\right\Vert _{MK}$%
\textit{, is a compact, hence complete, metric space, its topology being
exactly }$\mathcal{T(}w^{\ast },K^{n},a\mathcal{)}$ \textit{(in spite of the
fact that the normed space} $(cabv(K^{n}),\left\Vert .\right\Vert _{MK})$ 
\textit{is not complete if }$T$\textit{\ is infinite).}

\bigskip

\textbf{Remark}. The "basis" of\ Theorem 13 is Theorem 12 which asserts the
coincidence of convergent sequences in $\mathcal{T(}MK,K^{n},a\mathcal{)}$
and $\mathcal{T(}w^{\ast },K^{n},a\mathcal{)}$. This coincidence is no
longer valid for general $X$ instead of $K^{n}$ as we shall see later.

\bigskip

\textbf{The modified Monge-Kantorovich norm}

\bigskip

In this subparagraph, we shall be concerned with the so called "modified
Monge-Kantorovich norm", which can be defined only on a subspace of $cabv(X)$%
. This new norm is strongly related to the Monge-Kantorovich norm and
generates a most important distance (which generalizes classical
Kantorovich-Rubinstein metric on the space of probabilities, see e.g [6]) on
some distinguished subsets of $cabv(X)$.

For any $v\in X$, let us define 
\begin{equation*}
cabv(X,v)=\{\mu \in cabv(X)\mid \mu (T)=v\}\text{.}
\end{equation*}

Clearly $\delta _{t}v\in cabv(X,v)$ for any $t\in T$. It is seen that $%
cabv(X,0)$ is a vector subspace of $cabv(X)$. For any $\emptyset \neq
A\subset cabv(X,v)$ one has $A-A\overset{def}{=}\{\mu -\nu \mid \mu ,\nu \in
A\}\subset cabv(X,0)$.

\bigskip

\textbf{Lemma 14}. \textit{For any }$v\in X$\textit{, the set }$cabv(X,v)$%
\textit{\ is weak}$^{\ast }$\textit{\ closed in }$cabv(X)$\textit{.}

\textbf{Proof}. Take arbitrarily $x\in X$ and let us define the constant
function $\varphi _{x}:T\rightarrow X$, given via $\varphi _{x}(t)=x$, for
any $t\in T$.

Now take an arbitrary adherent point $\mu \in cabv(X)$ for $cabv(X,v)$.
Hence one can find $(\mu _{\delta })_{\delta }$ net $cabv(X,v)$ such that $%
\mu _{\delta }\underset{\delta }{\rightarrow }\mu $ in the topology $%
\mathcal{T(}w^{\ast },X)$, i.e. $\int fd\mu _{\delta }\underset{\delta }{%
\rightarrow }\int fd\mu $ for any $f\in C(X)$. Consequently $\int \varphi
_{x}d\mu _{\delta }\underset{\delta }{\rightarrow }\int \varphi _{x}d\mu $
for any $x\in X$, which means $(x\mid \mu _{\delta }(T))=(x\mid v)\underset{%
\delta }{\rightarrow }(x\mid \mu (T))$. So $(x\mid \mu (T))=(x\mid v)$ for
any $x\in X$, hence $\mu (T)=v$ and $\mu \in cabv(X,v)$. $\square $

\bigskip

Define%
\begin{equation*}
L_{1}(X)=\{f\in L(X)\mid \left\Vert f\right\Vert _{L}\leq 1\}
\end{equation*}%
and clearly $BL_{1}(X)\subset L_{1}(X)$.

For any $\mu \in cabv(X,0)$, let us define%
\begin{equation*}
\left\Vert \mu \right\Vert _{MK}^{\ast }\overset{def}{=}\sup \{\left\vert
\int fd\mu \right\vert \mid f\in L_{1}(X)\}\text{.}
\end{equation*}

\bigskip

\textbf{Theorem 15}. \textit{For any }$\mu \in cabv(X,0)$\textit{, one has}%
\begin{equation*}
\left\Vert \mu \right\Vert _{MK}\leq \left\Vert \mu \right\Vert _{MK}^{\ast
}\leq \left\Vert \mu \right\Vert diam(T)\text{.}
\end{equation*}

\textbf{Proof}. The first inequality is given by the inclusion $%
BL_{1}(X)\subset L_{1}(X)$.

To prove the second inequality, let us take arbitrarily $f\in L_{1}(X)$. For
any $t_{0}\in T$, one has%
\begin{equation*}
\left\vert \int fd\mu \right\vert =\left\vert \int (f-f(t_{0}))d\mu +\int
f(t_{0})d\mu \right\vert =
\end{equation*}%
\begin{equation*}
=\left\vert \int (f-f(t_{0}))d\mu +(f(t_{0})\mid \mu (T))\right\vert
=\left\vert \int (f-f(t_{0}))d\mu \right\vert \leq
\end{equation*}%
\begin{equation*}
\leq \left\Vert f-f(t_{0})\right\Vert \left\Vert \mu \right\Vert \text{.}
\end{equation*}

Because $\left\Vert f(t)-f(t_{0})\right\Vert \leq diam(T)$, for any $t\in T$%
, one has $\left\Vert f-f(t_{0})\right\Vert \leq diam(T)$, so $\left\vert
\int fd\mu \right\vert \leq \left\Vert \mu \right\Vert diam(T)$. $\square $

\bigskip

\textbf{Remark}. According to the definition, one has for any $\mu \in
cabv(X,0)$ and any $f\in L(X)$:%
\begin{equation}
\left\vert \int fd\mu \right\vert \leq \left\Vert \mu \right\Vert
_{MK}^{\ast }\left\Vert f\right\Vert _{L}\text{.}  \tag{10}
\end{equation}

Indeed, in case $\left\Vert f\right\Vert _{L}=0$, i.e. $f\equiv v\in X$ ($f$
is constant), one has $\int fd\mu =(v\mid \mu (T))=0$. In case $\left\Vert
f\right\Vert _{L}>0$, take $g=\frac{1}{\left\Vert f\right\Vert _{L}}f$ and $%
g\in L_{1}(X)$, hence $\left\vert \int gd\mu \right\vert \leq \left\Vert \mu
\right\Vert _{MK}^{\ast }$ a.s.o.

\bigskip

\textbf{Theorem 16}. \textit{The function} $p:cabv(X,0)\rightarrow \mathbb{R}%
_{+}$ \textit{given via }$p(\mu )=\left\Vert \mu \right\Vert _{MK}^{\ast }$%
\textit{\ is a norm on }$cabv(X,0)$\textit{.}

\textbf{Proof}. Using Theorem 15, one can see that $p$ takes finite values
and that $p(\mu )=0$ if and only if $\mu =0$. The fact that $p$ is a
seminorm is obvious. $\square $

\bigskip

\textbf{Definition 17}. The norm $\left\Vert .\right\Vert _{MK}^{\ast }$
defined above on $cabv(X,0)$ is called \textit{the modified
Monge-Kantorovich norm}.

\bigskip

\textbf{Theorem 18}. \textit{The norms }$\left\Vert .\right\Vert _{MK}$%
\textit{\ and }$\left\Vert .\right\Vert _{MK}^{\ast }$ \textit{are
equivalent on }$cabv(X,0)$\textit{. More precisely, for any }$\mu \in
cabv(X,0)$\textit{, one has}%
\begin{equation*}
\left\Vert \mu \right\Vert _{MK}\leq \left\Vert \mu \right\Vert _{MK}^{\ast
}\leq \left\Vert \mu \right\Vert _{MK}(diam(T)+1)
\end{equation*}%
\textit{and}%
\begin{equation*}
\left\Vert \mu \right\Vert _{MK}\leq \left\Vert \mu \right\Vert _{MK}^{\ast
}\leq \left\Vert \mu \right\Vert diam(T)
\end{equation*}

\textbf{Proof}. Take $\mu \in cabv(X,0)$. It remains to be proved that%
\begin{equation*}
\left\Vert \mu \right\Vert _{MK}^{\ast }\leq \left\Vert \mu \right\Vert
_{MK}(diam(T)+1)\text{.}
\end{equation*}

For arbitrary $f\in L_{1}(X)$ and $t_{0}\in T$, define $h:T\rightarrow K$
via $h(t)=f(t)-f(t_{0})$. Then $\left\Vert h\right\Vert _{L}=\left\Vert
f\right\Vert _{L}\leq 1$ and (obviously) $\left\Vert h\right\Vert \leq
diam(T)$. Consequently $\left\Vert h\right\Vert _{BL}\leq diam(T)+1$.
Because $\mu (T)=0$, one has $\int fd\mu =\int hd\mu $, hence (see (2))%
\begin{equation*}
\left\vert \int fd\mu \right\vert =\left\vert \int hd\mu \right\vert \leq
\left\Vert \mu \right\Vert _{MK}\left\Vert h\right\Vert _{BL}\leq \left\Vert
\mu \right\Vert _{MK}(diam(T)+1)
\end{equation*}%
and $f$ is arbitrary. $\square $

\bigskip

Theorem 18 says that the topology $\mathcal{T(}MK^{\ast },X)$ generated by $%
\left\Vert .\right\Vert _{MK}^{\ast }$\ on $cabv(X,0)$ coincides with he
topology induced by $\mathcal{T(}MK,X)$ on $cabv(X,0)$. For a sequence $(\mu
_{n})_{n\in \mathbb{N}}\subset cabv(X,0)$ and for $\mu \in cabv(X,0)$, one
has the equivalence: $\mu _{n}\underset{n}{\overset{\text{MK}^{\ast }}{%
\rightarrow }}\mu $ if and only if $\mu _{n}\underset{n}{\overset{\text{MK}}{%
\rightarrow }}\mu $.

\bigskip

\textbf{Theorem 19}. \textit{Let }$a$\textit{\ and }$b$\textit{\ be in }$T$%
\textit{, }$a\neq b$\textit{\ and }$x\in X$\textit{, }$\left\Vert
x\right\Vert =1$\textit{. Then }$\delta _{a}x-\delta _{b}x\in cabv(X,0)$%
\textit{\ and }%
\begin{equation*}
\left\Vert \delta _{a}x-\delta _{b}x\right\Vert _{MK}^{\ast }=d(a,b)\text{.}
\end{equation*}

\textbf{Proof}. For any $f\in L_{1}(X)$, writing $\mu =\delta _{a}x-\delta
_{b}x$, one has:%
\begin{equation*}
\left\vert \int fd\mu \right\vert =\left\vert (f(a)\mid x)-(f(b)\mid
x)\right\vert \leq \left\Vert f(a)-f(b)\right\Vert \left\Vert x\right\Vert =
\end{equation*}%
\begin{equation*}
=\left\Vert f(a)-f(b)\right\Vert \leq d(a,b)\text{,}
\end{equation*}%
hence%
\begin{equation*}
\left\Vert \mu \right\Vert _{MK}^{\ast }\leq d(a,b)\text{.}
\end{equation*}

Define $f:T\rightarrow X$, via $f(t)=d(t,a)x$. Then $f\in L_{1}(X)$,
because, for $u$ and $v$ in $T$, one has:%
\begin{equation*}
\left\Vert f(u)-f(v)\right\Vert =\left\Vert (d(u,a)-d(v,a))x\right\Vert =
\end{equation*}%
\begin{equation*}
=\left\vert d(u,a)-d(v,a)\right\vert \leq d(u,v)\text{.}
\end{equation*}%
Consequently%
\begin{equation*}
\left\vert \int fd\mu \right\vert =\left\vert (f(a)-f(b)\mid x)\right\vert =
\end{equation*}%
\begin{equation*}
\left\vert (-d(b,a)x\mid x)\right\vert =d(a,b)\leq \left\Vert \mu
\right\Vert _{MK}^{\ast }\text{. }\square
\end{equation*}

\bigskip

Notice that, if $a>0$ and $v\in X$ are such that $\left\Vert v\right\Vert
\leq a$, then $B_{a}(X,v)\overset{def}{=}B_{a}(X)\cap cabv(X,v)\neq
\emptyset $ (because $\left\Vert \delta _{t}v\right\Vert =\left\Vert
v\right\Vert \leq a$ for any $t\in T$).

On a non empty set $A\subset cabv(X)$, one can consider the following
distances:

- \textit{The} \textit{variational distance} given via $d_{\left\Vert
.\right\Vert }(\mu ,\nu )=\left\Vert \mu -\nu \right\Vert $.

- \textit{The} \textit{Monge-Kantorovich distance} given via $d_{MK}(\mu
,\nu )=\left\Vert \mu -\nu \right\Vert _{MK}$.

- Assuming that $A-A\subset cabv(X,0)$, \textit{the} \textit{modified
Monge-Kantorovich distance} given via $d_{MK}^{\ast }(\mu ,\nu )=\left\Vert
\mu -\nu \right\Vert _{MK}^{\ast }$.

We shall mainly work in the particular case when $A=B_{a}(X,v)$ with $%
\left\Vert v\right\Vert \leq a$. On such $B_{a}(X,v)$ the last two distances
are equivalent (Theorem 17): for $\mu $ and $\nu $ in $B_{a}(X,v)$, one has%
\begin{equation*}
d_{MK}(\mu ,\nu )\leq d_{MK}^{\ast }(\mu ,\nu )\leq d_{MK}(\mu ,\nu
)(diam(T)+1)\text{.}
\end{equation*}

In the next paragraph, we shall consider on such a set $A$ another distance
(namely the Hanin distance).

Before passing further, it is our duty to lay stress upon the fact that,
maybe, a more honest name for the (modified) Monge-Kantorovich distance
would have been Kantorovich-Rubinstein distance or Lipschitz distance.

\bigskip

\textbf{Theorem 20}. \textit{Let }$a>0$\textit{\ and }$v\in X$\textit{\ be} 
\textit{such that }$\left\Vert v\right\Vert \leq a$.

1.\textit{\ The (non empty) set} $B_{a}(X,v)$ \textit{is weak}$^{\ast }$%
\textit{\ closed in }$B_{a}(X)$\textit{, hence }$B_{a}(X,v)$\textit{\ is weak%
}$^{\ast }$\textit{\ compact. On }$B_{a}(X,v)$\textit{, one has the
equivalent metrics }$d_{MK}$\textit{\ and }$d_{MK}^{\ast }$\textit{.}

2.\textit{\ For any} $n\in \mathbb{N}$ \textit{(working for} $X=K^{n}$%
\textit{), one has the supplementary result that} $B_{a}(K^{n},v)$\textit{,
equipped with one of the equivalent metrics }$d_{MK}$\textit{\ and }$%
d_{MK}^{\ast }$\textit{\ is a compact, hence complete, metric space (its
topology being equal to the weak}$^{\ast }$\textit{\ topology on }$%
B_{a}(K^{n},v)$\textit{).}

3. \textit{In the particular case when }$K=\mathbb{R}$\textit{, }$n=1$%
\textit{\ and }$v\geq 0$\textit{, one can consider the set} $B_{a}^{+}(%
\mathbb{R},v)=B_{a}(\mathbb{R},v)\cap cabv^{+}(\mathbb{R})$, \textit{where }$%
cabv^{+}(\mathbb{R})=\{\mu \in cabv(\mathbb{R})\mid \mu \geq 0\}$\textit{.
Then }$B_{a}^{+}(\mathbb{R},v)$\textit{, equipped with one of the equivalent
metrics }$d_{MK}$\textit{\ and }$d_{MK}^{\ast }$\textit{\ is a compact,
hence complete, metric space (its topology being equal to the weak}$^{\ast }$%
\textit{\ topology on }$B_{a}^{+}(\mathbb{R},v)$\textit{).}

\textit{For} $a=v=1$, $B_{1}^{+}(\mathbb{R},1)$ \textit{is exactly the set
of all probabilities on} $\mathcal{B}$.

\textbf{Proof}. 1. We have $B_{a}(X,v)=B_{a}(X)\cap cabv(X,v)$. Because $%
B_{a}(X)$ is weak$^{\ast }$ compact, the result follows from Lemma 14.

2. The weak$^{\ast }$ topology of $B_{a}(K^{n})$ coincides with the topology
generated by the Monge-Kantorovich distance $d_{MK}$ (Theorem 13). Hence, $%
B_{a}(K^{n},v)$, being weak$^{\ast }$ closed in $B_{a}(K^{n})$ which is weak$%
^{\ast }$ compact, will be also weak$^{\ast }$ compact. Therefore $%
B_{a}(K^{n},v)$ is a compact subset of $B_{a}(K^{n})$, considering on $%
B_{a}(K^{n})$ the topology generated by the Monge-Kantorovich distance.
Because the Monge-Kantorovich distance and the modified Monge-Kantorovich
distance are equivalent on $B_{a}(K^{n},v)$, it follows that $B_{a}(K^{n},v)$
equipped either with the Monge-Kantorovich distance, or with the modified
Monge-Kantorovich distance is a compact, hence complete, metric space.

3. In the particular case $K=\mathbb{R}$, $n=1$, $v\geq 0$, all it remains
to be proved is the fact that $cabv^{+}(\mathbb{R})$ is weak$^{\ast }$
closed. To this end, let $\mu \in cabv(\mathbb{R})$ be such that there
exists $(\mu _{\delta })_{\delta }$ net $cabv^{+}(\mathbb{R})$ with the
property that $\mu _{\delta }\underset{\delta }{\rightarrow }\mu $ in $%
\mathcal{T(}w^{\ast },\mathbb{R}\mathcal{)}$. This implies that for any $%
f\in C(\mathbb{R})$, $f\geq 0$, one has $\int fd\mu _{\delta }\underset{%
\delta }{\rightarrow }\int fd\mu $. Because $\int fd\mu _{\delta }\geq 0$
for any $\delta $, it follows that $\int fd\mu \geq 0$. We succeeded in
proving that for any $0\leq f\in C(\mathbb{R})$, one has $\int fd\mu \geq 0$%
. So the functional $x_{\mu }^{^{\prime }}\in C(\mathbb{R})^{^{\prime }}$,
given via $x_{\mu }^{^{\prime }}(f)=\int fd\mu $ is positive. The
Riesz-Kakutani theorem says that this is equivalent to the fact that $\mu $
is positive, i.e. $\mu \in cabv^{+}(\mathbb{R})$. $\square $

\bigskip

\textbf{The Hanin norm}

\bigskip

The problem with the modified Monge-Kantorovich norm is the fact that it
cannot be defined on the whole space $cabv(X)$. To be more precise, using
the notation from Theorem 16, if one tries to extend $p$ beyond $cabv(X,0)$,
one can obtain infinite values for the extension, as the following result
shows.

\bigskip

\textbf{Proposition 21}. \textit{Define} $p:cabv(X)\rightarrow \overline{%
\mathbb{R}_{+}}$ \textit{via} $p(\mu )=\sup \{\left\vert \int fd\mu
\right\vert \mid f\in L_{1}(X)\}$. \textit{Then }$p$\textit{\ is an extended
seminorm (i.e. }$p(\mu +\nu )\leq p(\mu )+p(\nu )$\textit{\ and }$p(\alpha
\mu )=\left\vert \alpha \right\vert p(\mu )$\textit{\ (with convention }$%
0\cdot \infty =0$\textit{) for any }$\mu ,\nu \in cabv(X)$\textit{\ and any }%
$\alpha \in K$\textit{.}

\textit{We have the equivalence (for }$\mu \in cabv(X)$\textit{): }$p(\mu
)<\infty $\textit{\ }$\Leftrightarrow $\textit{\ }$\mu \in cabv(X,0)$)%
\textit{.}

\textbf{Proof}. The only fact which must be proved is the implication $%
\Rightarrow $ in the enunciation. Let $\mu \in cabv(X)$ with $p(\mu )<\infty 
$. Accepting $\mu (T)\neq 0$, we shall arrive at a contradiction.

Indeed, let $x\in X$ with $\left\Vert x\right\Vert =1$ and such that $(x\mid
\mu (T))=\left\Vert \mu (T)\right\Vert >0$. Then, for any $n\in \mathbb{N}$,
the function $f_{n}\in C(X)$ given via $f_{n}(t)=nx$ for any $t\in T$ is
constant, so $\left\Vert f_{n}\right\Vert _{L}=0$ and $\int f_{n}d\mu
=n\left\Vert \mu (T)\right\Vert \underset{n}{\rightarrow }\infty $. Hence,
because all $f_{n}\in L_{1}(X)$, it follows that $p(\mu )=\infty $ which is
a contradiction. $\square $

\bigskip

L. Hanin ([11] and [10]) succeeded in "extending" the modified
Monge-Kantorovich norm $\left\Vert .\right\Vert _{MK}^{\ast }$ from $cabv(%
\mathbb{R},0)$ to the whole $cabv(\mathbb{R)}$ (the "extension" is
equivalent to $\left\Vert .\right\Vert _{MK}^{\ast }$ on $cabv(\mathbb{R},0)$%
). Using approximately the same line of reasoning, we shall "extend" $%
\left\Vert .\right\Vert _{MK}^{\ast }$ from $cabv(X,0)$ to a norm $%
\left\Vert .\right\Vert _{H}$ on $cabv(X)$, for an arbitrary Hilbert space $%
X $.

Define for any $\mu \in cabv(X)$:%
\begin{equation*}
\left\Vert \mu \right\Vert _{H}\overset{def}{=}\inf \{\left\Vert \nu
\right\Vert _{MK}^{\ast }+\left\Vert \mu -\nu \right\Vert \mid \nu \in
cabv(X,0)\}\text{,}
\end{equation*}%
thus obtaining the map $\left\Vert .\right\Vert _{H}:cabv(X)\rightarrow 
\mathbb{R}_{+}$.

Taking $\nu =0$, we get $\left\Vert \mu \right\Vert _{H}\leq \left\Vert \mu
\right\Vert $ for any $\mu \in cabv(X)$. If $\mu \in cabv(X,0)$, taking $\nu
=\mu $, we obtain $\left\Vert \mu \right\Vert _{H}\leq \left\Vert \mu
\right\Vert _{MK}^{\ast }$.

\bigskip

\textbf{Proposition 22}. \textit{For any }$\mu \in cabv(X,0)$\textit{, one
has}%
\begin{equation*}
\left\Vert \mu \right\Vert _{MK}^{\ast }\leq \max (diam(T),1)\left\Vert \mu
\right\Vert _{H}\text{.}
\end{equation*}

\textbf{Proof}.\ For any $\nu \in cabv(X,0)$, using Theorem 18 we have:%
\begin{equation*}
\left\Vert \mu \right\Vert _{MK}^{\ast }=\left\Vert \nu +\mu -\nu
\right\Vert _{MK}^{\ast }\leq \left\Vert \nu \right\Vert _{MK}^{\ast
}+\left\Vert \mu -\nu \right\Vert _{MK}^{\ast }\leq
\end{equation*}%
\begin{equation*}
\leq \left\Vert \nu \right\Vert _{MK}^{\ast }+\left\Vert \mu -\nu
\right\Vert diam(T)=
\end{equation*}%
\begin{equation*}
=\left\Vert \nu \right\Vert _{MK}^{\ast }+\left\Vert \mu -\nu \right\Vert
+\left\Vert \mu -\nu \right\Vert (diam(T)-1)\text{.}
\end{equation*}

In case $diam(T)\leq 1$, we get, for any $\nu \in cabv(X,0)$: $\left\Vert
\mu \right\Vert _{MK}^{\ast }\leq \left\Vert \nu \right\Vert _{MK}^{\ast
}+\left\Vert \mu -\nu \right\Vert $ and passing to infimum we get $%
\left\Vert \mu \right\Vert _{MK}^{\ast }\leq \left\Vert \mu \right\Vert _{H}$%
.

In case $diam(T)>1$, we get, for any $\nu \in cabv(X,0)$:%
\begin{equation*}
\left\Vert \mu \right\Vert _{MK}^{\ast }\leq \left\Vert \nu \right\Vert
_{MK}^{\ast }+\left\Vert \mu -\nu \right\Vert +(\left\Vert \nu \right\Vert
_{MK}^{\ast }+\left\Vert \mu -\nu \right\Vert )(diam(T)-1)=
\end{equation*}%
\begin{equation*}
=(\left\Vert \nu \right\Vert _{MK}^{\ast }+\left\Vert \mu -\nu \right\Vert
)(1+diam(T)-1)
\end{equation*}%
and passing to infimum, we get $\left\Vert \mu \right\Vert _{MK}^{\ast }\leq
\left\Vert \mu \right\Vert _{H}diam(T)$. $\square $

\bigskip

Before passing further, let us notice that, for any $f\in L(X)$ and any $\mu
\in cabv(X)$, one has%
\begin{equation}
\left\vert \int fd\mu \right\vert \leq \left\Vert \mu \right\Vert
_{H}\left\Vert f\right\Vert _{BL}\text{.}  \tag{11}
\end{equation}

Indeed, for any $\nu \in cabv(X,0)$ we have:%
\begin{equation*}
\left\vert \int fd\mu \right\vert =\left\vert \int fd(\mu -\nu )+\int fd\nu
\right\vert \leq \left\vert \int fd\nu \right\vert +\left\vert \int fd(\mu
-\nu )\right\vert \leq
\end{equation*}%
\begin{equation*}
\leq \left\Vert \nu \right\Vert _{MK}^{\ast }\left\Vert f\right\Vert
_{L}+\left\Vert \mu -\nu \right\Vert \left\Vert f\right\Vert \leq
\end{equation*}%
\begin{equation*}
\leq \left\Vert f\right\Vert _{BL}(\left\Vert \nu \right\Vert _{MK}^{\ast
}+\left\Vert \mu -\nu \right\Vert )\text{,}
\end{equation*}%
with $(10)$. Due to the arbitrariness of $\nu $, $(11)$ is proved.

\bigskip

\textbf{Theorem 23}. 1. \textit{The functional} $\left\Vert .\right\Vert
_{H}:cabv(X)\rightarrow \mathbb{R}_{+}$ \textit{is a norm on }$cabv(X)$%
\textit{\ which generates a topology weaker than the variational topology
generated by }$\left\Vert .\right\Vert $\textit{: }$\left\Vert \mu
\right\Vert _{H}\leq \left\Vert \mu \right\Vert $ \textit{for any }$\mu \in
cabv(X)$.

2. \textit{On} $cabv(X,0)$\textit{, the modified Monge-Kantorovich norm }$%
\left\Vert .\right\Vert _{MK}^{\ast }$ \textit{and the restriction of} $%
\left\Vert .\right\Vert _{H}$ \textit{are equivalent: }%
\begin{equation*}
\left\Vert \mu \right\Vert _{H}\leq \left\Vert \mu \right\Vert _{MK}^{\ast
}\leq \max (diam(T),1)\left\Vert \mu \right\Vert _{H}
\end{equation*}%
\textit{\ (hence }$\left\Vert \mu \right\Vert _{H}=\left\Vert \mu
\right\Vert _{MK}^{\ast }$\textit{, if }$diam(T)\leq 1$\textit{) for any }$%
\mu \in cabv(X,0)$\textit{.}

\textbf{Proof}. It remains to be proved that $\left\Vert .\right\Vert _{H}$
is a norm on $cabv(X,0)$. First we prove that $\left\Vert .\right\Vert _{H}$
is a seminorm.

Because $cabv(X,0)=cabv(X,0)+cabv(X,0)$, we have for any $\mu _{1}$ and $\mu
_{2}$ in $cabv(X)$:%
\begin{equation*}
\left\Vert \mu _{1}+\mu _{2}\right\Vert _{H}=\inf \{\left\Vert \nu _{1}+\nu
_{2}\right\Vert _{MK}^{\ast }+\left\Vert \mu _{1}+\mu _{2}-\nu _{1}-\nu
_{2}\right\Vert \mid \nu _{1},\nu _{2}\in cabv(X,0)\}\text{.}
\end{equation*}

Because%
\begin{equation*}
\left\Vert \nu _{1}+\nu _{2}\right\Vert _{MK}^{\ast }+\left\Vert \mu
_{1}+\mu _{2}-\nu _{1}-\nu _{2}\right\Vert \leq
\end{equation*}%
\begin{equation*}
\leq \left\Vert \nu _{1}\right\Vert _{MK}^{\ast }++\left\Vert \nu
_{2}\right\Vert _{MK}^{\ast }+\left\Vert \mu _{1}-\nu _{1}\right\Vert
+\left\Vert \mu _{2}-\nu _{2}\right\Vert \text{,}
\end{equation*}%
for any $\nu _{1},\nu _{2}\in cabv(X,0)$, we pass to infimum obtaining $%
\left\Vert \mu _{1}+\mu _{2}\right\Vert _{H}\leq \left\Vert \mu
_{1}\right\Vert _{H}+\left\Vert \mu _{2}\right\Vert _{H}$.

For $\alpha \in K$ and $\mu \in cabv(X)$, one has $\left\Vert \alpha \mu
\right\Vert _{H}=\left\vert \alpha \right\vert \left\Vert \mu \right\Vert
_{H}$. This is obvious for $\alpha =0$. If $\alpha \neq 0$, using the
equality $\alpha cabv(X,0)=\{\alpha \mu \mid \mu \in cabv(X,0)\}=cabv(X,0)$
we have:%
\begin{equation*}
\left\Vert \alpha \mu \right\Vert _{H}=\inf \{\left\Vert \nu \right\Vert
_{MK}^{\ast }+\left\Vert \alpha \mu -\nu \right\Vert \mid \nu \in
cabv(X,0)\}=
\end{equation*}%
\begin{equation*}
=\inf \{\left\Vert \alpha \nu _{1}\right\Vert _{MK}^{\ast }+\left\Vert
\alpha \mu -\alpha \nu _{1}\right\Vert \mid \nu _{1}\in cabv(X,0)\}=
\end{equation*}%
\begin{equation*}
=\left\vert \alpha \right\vert \inf \{\left\Vert \nu _{1}\right\Vert
_{MK}^{\ast }+\left\Vert \mu -\nu _{1}\right\Vert \mid \nu _{1}\in
cabv(X,0)\}=\left\vert \alpha \right\vert \left\Vert \mu \right\Vert _{H}%
\text{.}
\end{equation*}

Finally, we show that, for $\mu \in cabv(X)$, one has the implication: $%
\left\Vert \mu \right\Vert _{H}=0\Rightarrow \mu =0$. Indeed, if $\left\Vert
\mu \right\Vert _{H}=0$, we have $\int fd\mu =0$ for any $f\in L(X)$, using $%
(11)$. This implies $\int fd\mu =0$ for any $f\in BL_{1}(X)$, hence $%
\left\Vert \mu \right\Vert _{MK}=0$, i.e. $\mu =0$. $\square $

\bigskip

\textbf{Corollary 24}. \textit{On} $cabv(X,0)$\textit{, the restrictions of} 
$\left\Vert .\right\Vert _{H}$ \textit{and }$\left\Vert .\right\Vert _{MK}$%
\textit{\ are equivalent: for any }$\mu \in cabv(X,0)$\textit{\ one has}%
\begin{equation*}
\frac{1}{\max (diam(T),1)}\left\Vert \mu \right\Vert _{MK}\leq \left\Vert
\mu \right\Vert _{H}\leq (diam(T)+1)\left\Vert \mu \right\Vert _{MK}\text{,}
\end{equation*}%
\textit{hence, if }$diam(T)\leq 1$\textit{, we hence}%
\begin{equation*}
\left\Vert \mu \right\Vert _{MK}\leq \left\Vert \mu \right\Vert _{H}\leq
2\left\Vert \mu \right\Vert _{MK}\text{.}
\end{equation*}

\textbf{Proof}. Using Theorems 18 and 23, we get for $\mu \in cabv(X,0)$:%
\begin{equation*}
\left\Vert \mu \right\Vert _{H}\leq \left\Vert \mu \right\Vert _{MK}^{\ast
}\leq (diam(T)+1)\left\Vert \mu \right\Vert _{MK}
\end{equation*}%
and%
\begin{equation*}
\left\Vert \mu \right\Vert _{H}\geq \frac{1}{\max (diam(T),1)}\left\Vert \mu
\right\Vert _{MK}^{\ast }\geq \frac{1}{\max (diam(T),1)}\left\Vert \mu
\right\Vert _{MK}\text{. }\square
\end{equation*}

\bigskip

\textbf{Definition 25}. The norm $\left\Vert .\right\Vert _{H}$ on $cabv(X)$
will be called \textit{the Hanin norm}.

\bigskip

As usual, we present the afferent notations. The topology on $cabv(X)$,
generated by $\left\Vert .\right\Vert _{H}$, will be called \textit{the} 
\textit{Hanin topology} and it will be denoted by $\mathcal{T(}H,X\mathcal{)}
$. This topology induces the topology $\mathcal{T(}H,X,a\mathcal{)}$ on $%
B_{a}(X)$. For a sequence $(\mu _{n})_{n}\subset cabv(X)$ and for $\mu \in
cabv(X)$, we write $\mu _{n}\overset{H}{\underset{n}{\rightarrow }}\mu $ to
denote the fact that $(\mu _{n})_{n}$ converges to $\mu $ in the topology $%
\mathcal{T(}H,X\mathcal{)}$. Finally, on any $\emptyset \neq A\subset
cabv(X) $, $\left\Vert .\right\Vert _{H}$ generates \textit{the Hanin metric}
$d_{H}$ given via $d_{H}(\mu ,\nu )=\left\Vert \mu -\nu \right\Vert _{H}$
for any $\mu ,\nu $ in $A$.

\bigskip

Theorem 23 is very important. Due to the equivalence of $\left\Vert
.\right\Vert _{H}$ and $\left\Vert .\right\Vert _{MK}^{\ast }$ on $cabv(X,0)$%
, it follows that that the metrics $d_{H}$ and $d_{MK}^{\ast }$ are
equivalent on any $B_{a}(X,v)$. Hence, in the enunciation of Theorem 20, one
can add $d_{H}$ to the previous equivalent metrics $d_{MK}$ and $%
d_{MK}^{\ast }$. Consequently we have:

\bigskip

\textbf{Corollary 26}. \textit{If }$a>0$\textit{\ and }$v\in X$\textit{\ are
such that }$\left\Vert v\right\Vert \leq a$\textit{, then the metrics }$%
d_{MK}$\textit{, }$d_{MK}^{\ast }$\textit{\ and }$d_{H}$\textit{\ are
equivalent on the non empty set }$B_{a}(X,v)$\textit{.}

\textit{Consequently, in the enunciation of Theorem 20, one can replace
(three times) the expression "the equivalent metrics }$d_{MK}$\textit{\ and }%
$d_{MK}^{\ast }$\textit{" with the expression "the equivalent metrics }$%
d_{MK}$\textit{, }$d_{MK}^{\ast }$\textit{\ and }$d_{H}$\textit{", thus
augmenting the enunciation.}

\bigskip

\textbf{Theorem 27}. \textit{Let }$a$\textit{\ and }$b$\textit{\ in }$T$%
\textit{, }$a\neq b$\textit{\ and }$x\in X$, $\left\Vert x\right\Vert =1$. 
\textit{Then}%
\begin{equation*}
\left\Vert \delta _{a}x\right\Vert _{H}=1
\end{equation*}%
\textit{and}%
\begin{equation*}
\frac{1}{\max (diam(T),1)}d(a,b)\leq \left\Vert \delta _{a}x-\delta
_{b}x\right\Vert _{H}\leq d(a,b)\text{.}
\end{equation*}

\textit{Hence, in case }$diam(T)\leq 1$\textit{, }$\left\Vert \delta
_{a}x-\delta _{b}x\right\Vert _{H}=d(a,b)$\textit{. For }$X=K$\textit{\ and }%
$x=1$\textit{, we have }%
\begin{equation*}
\frac{1}{\max (diam(T),1)}d(a,b)\leq \left\Vert \delta _{a}-\delta
_{b}\right\Vert _{H}\leq d(a,b)\text{.}
\end{equation*}

\textbf{Proof}. Again define $f=\varphi _{x}$, where $\varphi
_{x}:T\rightarrow X$, $\varphi _{x}(t)=x$ for any $t\in T$ and notice that $%
\left\Vert f\right\Vert _{BL}=\left\Vert f\right\Vert =1$.

For $\mu =\delta _{a}x$ use $(11)$ and get 
\begin{equation*}
\left\vert \int fd\mu \right\vert =(x\mid \mu (T))=(x\mid x)=1\leq
\left\Vert \mu \right\Vert _{H}\leq \left\Vert \mu \right\Vert =1.
\end{equation*}

Now, writing $\nu =\delta _{a}x-\delta _{b}x\in cabv(X,0)$ and using
Theorems 19 and 23, we get:%
\begin{equation*}
\left\Vert \nu \right\Vert _{H}\leq \left\Vert \nu \right\Vert _{MK}^{\ast
}=d(a,b)
\end{equation*}%
and%
\begin{equation*}
\left\Vert \nu \right\Vert _{H}\geq \frac{1}{\max (diam(T),1)}\left\Vert \nu
\right\Vert _{MK}^{\ast }=\frac{1}{\max (diam(T),1)}d(a,b)\text{. }\square
\end{equation*}

\bigskip

\textbf{Infinite dimensional extensions do not work}

\bigskip

In this subparagraph we present a counterexample showing that neither
Theorem 11, nor Theorem 12, can be extended for general Hilbert spaces
instead of $K^{n}$.

\bigskip

\textbf{Counterexample. }Our compact metric space\textbf{\ }$(T,d)$ will be
given as follows: $T=\{1,2\}$ with metric $d(i,j)=1$ if $i\neq j$ and $%
d(i,j)=0$ if $i=j$. Hence, on $T$ we have the discrete topology and the
Borel sets are $\mathcal{B}=\mathcal{P}(T)$. Our Hilbert space will be $%
l^{2} $.

Any function $f:T\rightarrow l^{2}$ is continuous, even Lipschitz, and
simple. We identify such a function $f$ giving $f(1)=(a_{1m})_{m}$ and $%
f(2)=(a_{2m})_{m}$. Clearly $\left\Vert f\right\Vert _{L}=\left\Vert
f(1)-f(2)\right\Vert $.

A measure $\mu \in cabv(l^{2})$ is identified giving $\mu
(\{1\})=(b_{1m})_{m}\in l^{2}$ and $\mu (\{2\})=(b_{2m})_{m}\in l^{2}$.
Hence, the total variation $\left\vert \mu \right\vert (T)=\left\Vert \mu
\right\Vert =\left\Vert \mu (\{1\})\right\Vert +\left\Vert \mu
(\{2\})\right\Vert $.

Using the previous notations, we have:

\begin{equation*}
\int fd\mu =(f(1)\mid \mu (\{1\}))+(f(2)\mid \mu (\{2\}))==\underset{m=1}{%
\overset{\infty }{\sum }}(a_{1m}\overline{b_{1m}}+a_{2m}\overline{b_{2m}})%
\text{.}
\end{equation*}

We shall be concerned with the case when $\mu \in cabv(l^{2},0)$, which
means $\mu (\{1\})+\mu (\{2\})=0$, i.e. $(b_{1m}+b_{2m})_{m}=0$, hence $%
b_{2m}=-b_{1m}$ for any $m$. In this case, we shall identify $\mu \equiv
b=(b_{m})_{m}\in l^{2}$, where $b_{1m}=b_{m}$ and $b_{2m}=-b_{m}$ for any $m$%
.

We have $\left\Vert \mu \right\Vert =\left\Vert b\right\Vert +\left\Vert
b\right\Vert =2\left\Vert b\right\Vert $ and $\int fd\mu =\underset{m=1}{%
\overset{\infty }{\sum }}(a_{1m}-a_{2m})\overline{b_{m}}$.

Writing $a=(a_{1m}-a_{2m})_{m}\in l^{2}$, we notice that $\int fd\mu =(a\mid
b)$. Notice also that $f\in L_{1}(l^{2})$ means $\left\Vert
f(1)-f(2)\right\Vert \leq 1$ i.e. $\left\Vert a\right\Vert \leq 1$.

The final preliminary fact is that for $b\equiv \mu \in cabv(l^{2},0)$, one
has $\left\Vert \mu \right\Vert _{MK}^{\ast }=\left\Vert b\right\Vert $.
Indeed $\left\Vert \mu \right\Vert _{MK}^{\ast }=\sup \{\left\vert \int
fd\mu \right\vert \mid \left\Vert f\right\Vert _{L}\leq 1\}$. For $f$
identified as above, we saw that $\left\Vert f\right\Vert _{L}\leq 1$ means $%
\left\Vert a\right\Vert \leq 1$, hence $\left\vert \int fd\mu \right\vert
=\left\vert (a\mid b)\right\vert \leq \left\Vert b\right\Vert $,
consequently $\left\Vert \mu \right\Vert _{MK}^{\ast }\leq \left\Vert
b\right\Vert $. On the other hand, let us take $a=(a_{m})_{m}\in l^{2}$ such
that $\left\Vert a\right\Vert =1$ and $(a\mid b)=\left\Vert b\right\Vert $.
Define $a_{1m}=a_{m}$ and $a_{2m}=0$, for every $m\in \mathbb{N}$. We got
the function $f:T\rightarrow l^{2}$ identified as above: $%
f(1)=(a_{1m})_{m}=(a_{m})_{m}$ and $f(2)=(a_{2m})_{m}=0$. Then, for this $f$
one has $\left\Vert f(1)-f(2)\right\Vert =\left\Vert a\right\Vert =1$ (hence 
$f\in L_{1}(l^{2})$) and $\int fd\mu =(a\mid b)=\left\Vert b\right\Vert $.
It follows that $\left\Vert b\right\Vert \leq \left\Vert \mu \right\Vert
_{MK}^{\ast }$ and the equality $\left\Vert \mu \right\Vert _{MK}^{\ast
}=\left\Vert b\right\Vert $ is proved.

Practically, we proved the existence of the linear and isometric isomorphism 
$(cabv(l^{2},0),\left\Vert .\right\Vert _{MK}^{\ast })\equiv l^{2}$, via $%
\mu \equiv b$ as above.

Let us consider an arbitrary bounded sequence $(b^{m})_{m}$ in $l^{2}$ and
let $a>0$ be such that $\left\Vert b^{m}\right\Vert \leq \frac{a}{2}$ for
any $m$. According to the previous isomorphism considerations, identify each 
$b^{m}\equiv \mu ^{m}\in cabv(l^{2},0)$. We saw that $\left\Vert \mu
^{m}\right\Vert =2\left\Vert b^{m}\right\Vert \leq a$, hence $\mu ^{m}\in
B_{a}(l^{2})$ for any $m$.

According to Theorem 18 and previous considerations we have the
equivalences: $\mu ^{m}\overset{\text{MK}}{\underset{m}{\rightarrow }}0$ if
and only if $\mu ^{m}\overset{\text{MK}^{\ast }}{\underset{m}{\rightarrow }}%
0 $ if and only $b^{m}\underset{m}{\rightarrow }0$, the last convergence
being in $l^{2}$.

On the other hand, we have the equivalences: $\mu ^{m}\overset{\text{w}%
^{\ast }}{\underset{m}{\rightarrow }}0$ if and only if\linebreak\ $(a\mid
b^{m})\underset{m}{\rightarrow }0$ for any $a\in l^{2}$ if and only if $b^{m}%
\overset{\text{w}^{\ast }}{\underset{m}{\rightarrow }}0$ (the last
convergence being weak convergence in $l^{2}$). We must prove the
equivalence: $\mu ^{m}\overset{\text{w}^{\ast }}{\underset{m}{\rightarrow }}%
0 $ $\Leftrightarrow (a\mid b^{m})\underset{m}{\rightarrow }0$ for any $a\in
l^{2}$. To prove $"\Rightarrow "$, take $a=(a_{n})_{n}\in l^{2}$ and define $%
f\in C(l^{2})$ identified via $f(1)=(a_{n})_{n}$ and $f(2)=0$. Then $(a\mid
b^{m})=\int fd\mu ^{m}\underset{m}{\rightarrow }0$. To prove $"\Leftarrow "$%
, take $f\in C(l^{2})$ identified via $f(1)=(a_{1n})_{n}$ and $%
f(2)=(a_{2n})_{n}$ and define $a=(a_{1n}-a_{2n})_{n}\in l^{2}$. Then $\int
fd\mu ^{m}=(a\mid b^{m})\underset{m}{\rightarrow }0$.

We arrived at the end of our discussion. Accept that the assertion in
Theorem 12 is valid for the separable Hilbert space $l^{2}$ instead of $%
K^{n} $ for any number $a>0$. This means to accept the following assertion
(taking $\mu ^{m}-\mu $ instead of $\mu ^{m}$ and $0$ instead of $\mu $ in
Theorem 12): For any number $a>0$ and for any sequence $(\mu
^{m})_{m}\subset B_{a}(l^{2},0)$, one has the equivalence: $\mu ^{m}\overset{%
\text{MK}^{\ast }}{\underset{m}{\rightarrow }}0$ if and only if $\mu ^{m}%
\overset{\text{w}^{\ast }}{\underset{m}{\rightarrow }}0$. This last
assertion "translated" in view of the previous considerations means: For any
number $a>0$ and for any bounded sequence $(b^{m})_{m}\subset l^{2}$ such
that $\left\Vert b^{m}\right\Vert \leq \frac{a}{2}$ for any $m$, one has the
equivalence: $b^{m}\underset{m}{\rightarrow }0$ (in $l^{2}$) if and only if $%
b^{m}\underset{m}{\rightarrow }0$ (weakly in $l^{2}$).

The last equivalence is clearly false. For a concrete example of failure,
one can take $a=\frac{2\pi }{\sqrt{6}}$ and the sequence $(b^{p})_{p}$ with $%
b^{p}=(b_{n}^{p})_{n}$, as follows:

For $p=1$: $b_{n}^{1}=\frac{1}{n}$, $n\in \mathbb{N}$.

For $p>1$: $b_{n}^{p}=\{%
\begin{array}{cc}
0\text{,} & \text{if }n<p \\ 
\frac{1}{n-p+1}\text{,} & \text{if }n\geq p%
\end{array}%
$.

Then $\left\Vert b^{p}\right\Vert =\frac{\pi }{\sqrt{6}}=\frac{a}{2}$ for
any $p$. Clearly the assertion $b^{p}\underset{p}{\rightarrow }0$ in $l^{2}$
is false.

On the other hand, let us take $a=(a_{n})_{n}\in l^{2}$. Defining, for any $%
p\in \mathbb{N}$, the new sequence $a(p)=(a_{n+p-1})_{n}$, one has $%
\left\Vert a(p)\right\Vert ^{2}=\underset{n=1}{\overset{\infty }{\sum }}%
\left\vert a_{n+p-1}\right\vert ^{2}\underset{p}{\rightarrow }0$. For any $p$%
, one can see that $(a\mid b^{p})=(a(p)\mid b^{1})$. It follows that $%
\left\vert (a\mid b^{p})\right\vert \leq \left\Vert a(p)\right\Vert
\left\Vert b^{1}\right\Vert \underset{p}{\rightarrow }0$ and this shows that 
$b^{p}\underset{p}{\rightarrow }0$ weakly in $l^{2}$.

\bigskip

\textbf{Remark}. The previous counterexample shows that the "extension" of
Theorem 12 for an arbitrary $X$ instead of $K^{n}$ is false. At the same
time, the counterexample shows that the "extension" of Theorem 11 for an
arbitrary $X$ instead of $K^{n}$ is false too. Indeed, accepting for
instance that $BL_{1}(l^{2})$ is relatively compact in $C(l^{2})$, we can
repeat the "proof" of Theorem 12 for $l^{2}$ instead of $K^{n}$ and arrive
at the false conclusion that weak$^{\ast }$ convergence and
Monge-Kantorovich convergence are the same in $cabv(l^{2})$.

\bigskip

\textbf{New metrics on }$T$ \textbf{generated by the previous norms on }$%
cabv(X,0)$

\bigskip

In this last subparagraph, we shall discuss the new metrics on $T$ which are
generated by the norms on $cabv(X,0)$ which have been introduced throughout
the paper.

To begin, let us choose an arbitrary $x\in X$ with $\left\Vert x\right\Vert
=1$ which will be fixed from now on (in the special case $X=K$, we take
canonically $x=1$). Recall that $\left\Vert \delta _{t}x\right\Vert
=\left\Vert \delta _{t}x\right\Vert _{MK}=\left\Vert \delta _{t}x\right\Vert
_{H}=1$ (relation $(3)$ and Theorem 27). We define the injective map $%
V:T\rightarrow cabv(X)$ via $V(t)=\delta _{t}x$. Then $\delta _{t}x-\delta
_{s}x\in cabv(X,0)$ and, for any norm $p$ on $cabv(X,0)$, one obtains the
metric $\rho _{p}$ on $T$ given via $\rho _{p}(t,s)=p(\delta _{t}x-\delta
_{s}x)$.

The fact that $\rho _{\left\Vert .\right\Vert }(t,s)=2$ for $t\neq s$ shows
that the metric $\rho _{\left\Vert .\right\Vert }$ generates the discrete
topology on $T\,$, being metrically insensitive (rigid) (see relation $(4)$).

To complete our discussion, we shall need.

\bigskip

\textbf{Lemma 28}. \textit{For any }$a\neq b$\textit{\ in }$T$\textit{,
there exists} $g\in L(\mathbb{R})$ \textit{having the following properties:}

\textit{i) }$g(a)=0$\textit{, }$g(b)=1$\textit{\ and }$0\leq g(t)\leq 1$%
\textit{\ for any }$t\in T$\textit{.}

\textit{ii) }$\left\Vert g\right\Vert =1$, $\left\Vert g\right\Vert _{L}\leq 
\frac{1}{d(a,b)}$\textit{, hence }$\left\Vert g\right\Vert _{BL}\leq \frac{%
1+d(a,b)}{d(a,b)}$\textit{.}

\textbf{Proof}.\textit{\ }Let us define\textit{\ }$g:T\rightarrow \mathbb{R}$%
, via\textit{\ }%
\begin{equation*}
g(t)=\frac{d(t,a)}{d(t,a)+d(t,b)}\text{.}
\end{equation*}

Then $i)$ and $\left\Vert g\right\Vert =1$ follow immediately. It remains to
be proved that $\left\Vert g\right\Vert _{L}\leq \frac{1}{d(a,b)}$.

Indeed, for any $x$ and $y$ in $T$, one has%
\begin{equation*}
\left\vert g(x)-g(y)\right\vert =\frac{\left\vert
d(x,a)d(y,b)-d(x,b)d(y,a)\right\vert }{(d(x,a)+d(x,b))(d(y,a)+d(y,b))}\leq
\end{equation*}%
\begin{equation*}
\leq \frac{\left\vert d(x,a)d(y,b)-d(x,b)d(y,a)\right\vert }{%
d(a,b)(d(y,a)+d(y,b))}=
\end{equation*}%
\begin{equation*}
=\frac{\left\vert
d(x,a)d(y,b)-d(y,a)d(y,b)+d(y,a)d(y,b)-d(y,a)d(x,b)\right\vert }{%
d(a,b)(d(y,a)+d(y,b))}\leq
\end{equation*}%
\begin{equation*}
\leq \frac{d(y,b)\left\vert d(x,a)-d(y,a)\right\vert +d(y,a)\left\vert
d(x,b)-d(y,b)\right\vert }{d(a,b)(d(y,a)+d(y,b))}\leq
\end{equation*}%
\begin{equation*}
\leq \frac{d(y,b)d(x,y)+d(y,a)d(x,y)}{d(a,b)(d(y,a)+d(y,b))}\leq \frac{d(x,y)%
}{d(a,b)}\text{. }\square
\end{equation*}

\bigskip

\textbf{Theorem 29}. \textit{The metrics }$d$\textit{, }$\rho _{\left\Vert
.\right\Vert _{MK}}$\textit{, }$\rho _{\left\Vert .\right\Vert _{MK}^{\ast
}} $ \textit{and }$\rho _{\left\Vert .\right\Vert _{H}}$\textit{\ are
Lipschitz equivalent. Namely, one has}%
\begin{equation}
\frac{1}{1+diam(T)}d\leq \rho _{\left\Vert .\right\Vert _{MK}}\leq d\text{,}
\tag{12}
\end{equation}

\begin{equation}
\rho _{\left\Vert .\right\Vert _{MK}^{\ast }}=d  \tag{13}
\end{equation}%
\textit{and}%
\begin{equation}
\frac{1}{\max (diam(T),1)}d\leq \rho _{\left\Vert .\right\Vert _{H}}\leq d%
\text{.}  \tag{14}
\end{equation}

\textit{The metric }$\rho _{\left\Vert .\right\Vert }$\textit{\ generates
the discrete topology on }$T$\textit{, being equivalent to }$d$\textit{\ if
and only if }$T$\textit{\ is finite.}

\textbf{Proof}. Relation $(13)$ is proved in Theorem 19 and relation $(14)$
is proved in Theorem 27. Because $\rho _{\left\Vert .\right\Vert _{MK}}\leq
d $ (see relation $(5)$), all it remains to be proved is the fact
(completely proving $(12)$) that 
\begin{equation}
\frac{1}{1+diam(T)}d\leq \rho _{\left\Vert .\right\Vert _{MK}}\text{.} 
\tag{15}
\end{equation}

To this end, for the previous $x\in X$ with $\left\Vert x\right\Vert =1$,
define $f=gx$, where $g$ is the function from Lemma 28, for arbitrary $a\neq
b$ in $T$. Then, for any $\mu \in cabv(X)$ one has%
\begin{equation*}
\left\vert \int fd\mu \right\vert \leq \left\Vert \mu \right\Vert
_{MK}\left\Vert f\right\Vert _{BL}=\left\Vert \mu \right\Vert
_{MK}\left\Vert g\right\Vert _{BL}\leq
\end{equation*}%
\begin{equation*}
\leq \left\Vert \mu \right\Vert _{MK}\frac{1+d(a,b)}{d(a,b)}\leq \frac{%
1+diam(T)}{d(a,b)}\left\Vert \mu \right\Vert _{MK}
\end{equation*}%
and this implies 
\begin{equation*}
\left\Vert \mu \right\Vert _{MK}\geq \frac{d(a,b)}{1+diam(T)}\left\vert \int
fd\mu \right\vert \text{.}
\end{equation*}

Taking $\mu =\delta _{a}x-\delta _{b}x$, one has $\int fd\mu
=((f(a)-f(b))\mid x)=-(x\mid x)=-1$, hence $\left\Vert \mu \right\Vert
_{MK}=\rho _{\left\Vert .\right\Vert _{MK}}(a,b)\geq \frac{d(a,b)}{1+diam(T)}
$ and $(15)$ is proved. $\square $

\bigskip

\textbf{Remarks}. \textbf{1}. The equivalence between the metrics $d$ and
the metrics $\rho _{\left\Vert .\right\Vert _{MK}}$\textit{, }$\rho
_{\left\Vert .\right\Vert _{MK}^{\ast }}$ and\textit{\ }$\rho _{\left\Vert
.\right\Vert _{H}}$ can be obtained immediately observing the fact that $%
(T,d)$ is a metric compact space and the other metrics generate weaker
topologies. Of course, the result of Theorem 29 is stronger.

\textbf{2}. Because (see Lemma 28) one has, for any $a\neq b$ in $T$, the
inequality $\left\Vert g\right\Vert _{BL}\leq \frac{1+d(a,b)}{d(a,b)}$, the
proof of Theorem 29 shows that, actually, for any $a$ and $b$ in $T$: $\frac{%
1}{1+diam(T)}d(a,b)\leq \frac{1}{1+d(a,b)}d(a,b)\leq \rho _{\left\Vert
.\right\Vert _{MK}}(a,b)\leq d(a,b)$ thus improving $(12)$.

\textbf{3}. In some cases, the evaluations given in Theorem 29 can be
improved. For instance, let us consider a number $A>0$, obtaining the
compact metric space $T=[0,A]$ with metric $d$ given by the natural metric
of $\mathbb{R}$. For any $a$ and $b$ in $T$, we use, for the computation of $%
\left\Vert \delta _{a}-\delta _{b}\right\Vert _{MK}=\rho _{\left\Vert
.\right\Vert _{MK}}(a,b)$, affine functions $f:T\rightarrow \mathbb{R}$ of
the form $f(x)=mx+n$ and, with some effort, we can conclude that $\rho
_{\left\Vert .\right\Vert _{MK}}(a,b)\geq \frac{2\left\vert a-b\right\vert }{%
2+A}$. This means that, in this case, one has%
\begin{equation}
\frac{2}{2+diam(T)}d(a,b)\leq \rho _{\left\Vert .\right\Vert _{MK}}(a,b)%
\text{.}  \tag{16}
\end{equation}

Relation $(16)$ improves relation $(15)$, i.e. improves relation $(12)$.

At the same time, taking $A$ arbitrary small, one can conclude from $(16)$
that, for any $\varepsilon \in (0,1)$, there exists a compact metric space $%
(T,d)$ such that $(1-\varepsilon )d\leq \rho _{\left\Vert .\right\Vert
_{MK}}\leq d$. Hence, it is natural to try to solve the following

\textit{Open Problem}: Find new evaluations, sharper than $(12)$, $(13)$ and 
$(14)$.

\bigskip

\textbf{References}

\bigskip

[1] P. Appel, M\'{e}moire sur les d\'{e}blais et les remblais des syst\`{e}%
mes continus ou discontinus. M\'{e}moirs pr\'{e}sent\'{e}es par divers
Savants \`{a} l'Acad\'{e}mie des Sciences de l'Institut de France, Paris 29
(1887), 1-208.

[2] P. Billingsley, Convergence of Probability Measures, second ed., John
Wiley \& Sons, New York, 1999.

[3] I. Chi\c{t}escu, Function Spaces, Ed. \c{S}t. Encicl. Bucure\c{s}ti,
1983 (in Romanian).

[4] J. Diestel, J.J. Uhl, Jr., Vector Measures, Math. Surveys, vol. 15,
Amer. Math. Soc., Providence, R.I., 1977.

[5] N. Dinculeanu, Vector Measures, VEB Deutscher Verlag der
Wissenschafthen, Berlin, 1966.

[6] R.M. Dudley, Real Analysis and Probability, second ed., Cambridge
University Press, Cambridge, 2002.

[7] J. Dugundji, Topology, fourth ed., Allyn and Bacon, Boston, 1968.

[8] N. Dunford, J.T. Schwartz, Linear Operators, Part I: General Theory,
Interscience Publishers, New York, 1957.

[9] G.B. Folland, Real Analysis. Modern Techniques and Their Applications,
second ed., John Wiley \& Sons, New York, 1999.

[10] L.G. Hanin, Kantorovich-Rubinstein norm and its application in the
theory of Lipschitz spaces, Proc. Amer. Math. Soc. 115 (1992) 345-352.

[11] L.G. Hanin, An extension of the Kantorovich norm. In Monge-Amp\`{e}re
Equation: Applications to geometry and optimization (Deerfield Beach, FL,
1997), vol. 226 of Contemp. Math., Amer. Math. Soc., Providence, RI (1999),
113-130.

[12] L.V. Kantorovich, On the translocation of masses, C.R. (Doklady) Acad.
Sci. URSS 37 (1942) 199-201.

[13] L.V. Kantorovich, G.P. Akilov, Functional Analysis, second ed., Per-%
\newline
gamon Press, New York, 1982.

[14] L.V. Kantorovich, G.S. Rubinstein, On a certain function space and some
extremal problems, C.R. (Doklady) Acad. Sci. URSS 115 (1957) 1058-1061 (in
Russian).

[15] L.V. Kantorovich, G.S. Rubinstein, On a space of completely additive
functions, Vestnik Leningrad Univ. 13 (1958) 52-59 (in Russian).

[16] J.L. Kelley, General topology, D. Van Nostrand Company, Toronto-New
York-London, 1955.

[17] J. Luke\v{s}, J. Maly, Measure and Integral, Matfyzpress, Publishing
House of the Faculty of Mathematics and Physics, Charles University, Prague,
1995.

[18] G. Monge, M\'{e}moire sur la th\'{e}orie des d\'{e}blais et des
remblais. In Histoire de l'Acad\'{e}mie Royale des Sciences de Paris, 1781,
666-704.

[19] K.R. Parthasarathy, Probability measures on metric spaces, Academic
Press, New York, 1967.

[20] H.H. Schaefer, Topological vector spaces, third printing corrected,
Springer, New York, 1971.

[21] C. Villani, Topics in optimal transportation, Graduate Studies in
Mathematics, vol. 58, Amer. Math. Soc., Providence, R.I., 2003.

\end{document}